\documentclass[pdflatex,sn-mathphys-num]{sn-jnl}
\usepackage{amsmath,amssymb,amsthm,mathtools}
\usepackage{needspace,xcolor}
\hypersetup{hidelinks,
  pdftitle={Vanishing Dissipation Limit to Planar Euler Shock for the Full Compressible Navier-Stokes Equations in R3},
  pdfauthor={Yujie Wang}}
\newcommand{\eps}{\varepsilon}
\newcommand{\R}{\mathbb R}
\newcommand{\dd}{\,\mathrm d}
\newcommand{\Div}{\operatorname{div}}
\newcommand{\norm}[1]{\left\lVert#1\right\rVert}
\newcommand{\abs}[1]{\left\lvert#1\right\rvert}
\theoremstyle{thmstyleone}
\newtheorem{theorem}{Theorem}[section]
\newtheorem{proposition}[theorem]{Proposition}
\newtheorem{corollary}[theorem]{Corollary}
\newtheorem{lemma}[theorem]{Lemma}
\theoremstyle{thmstyletwo}
\newtheorem{remark}[theorem]{Remark}
\numberwithin{equation}{section}
\allowdisplaybreaks[3]
\AtBeginDocument{%
  \everydisplay\expandafter{\the\everydisplay
    \abovedisplayskip=0pt
    \belowdisplayskip=0pt
    \abovedisplayshortskip=0pt
    \belowdisplayshortskip=0pt
  }%
}

\begin{document}

\title[Vanishing Dissipation Limit to Planar Euler Shock]{Vanishing Dissipation Limit to Planar Euler Shock
for the Full Compressible Navier--Stokes Equations in $\R^3$}

\author*[1]{\fnm{Yujie} \sur{Wang}}
\email{15316072697@163.com}

\affil[1]{\orgdiv{School of Mathematics, Statistics and Mechanics}, \orgname{Beijing University of Technology}, \orgaddress{\city{Beijing}, \postcode{100124}, \country{China}}}

\abstract{
\unboldmath
We prove the vanishing dissipation limit to a weak planar Euler shock
for the full compressible Navier--Stokes equations in $\R^3$.
The viscosity and heat conductivity are proportional to $\eps$,
with their ratio fixed. The conservative initial perturbation may
depend on all three spatial variables and is assumed small only in
the scaled $H^2\cap\dot H^{-1}$ norm. No $L^1$ bound or zero-mass
condition is required. For these data, we obtain a unique global
strong solution and an energy estimate uniform in time and in $\eps$.
On every cylinder with a measurable transverse section of finite
positive area, the $L^2$ error has an explicit leading term,
yielding the sharp rate $O(\sqrt\eps)$. The remainder has order
$\eps^{3/2}$ times the scaled initial norm; both estimates are
uniform for $t\ge0$.
If this norm tends to zero, convergence away from the shock front
is uniform in time, including at $t=0$. The entropy production
measures converge locally weakly to the measure carried by the
Euler shock. Low frequencies are controlled by an entropy-weighted
elliptic estimate. At second order, integration by parts in the
pressure and viscous heating terms leaves at most two derivatives
on the density, while third derivatives of velocity and temperature
enter only through the time-integrated dissipation.

}


\maketitle

\section{Introduction}

A viscous profile connects the end states of an Euler shock
through a smooth layer of thickness $O(\eps)$. For perturbations
depending on all three spatial variables, we establish stability
estimates uniform for $t\ge0$ as $\eps\to0$.

We consider the full compressible Navier--Stokes equations in three
space dimensions on $[0,\infty)\times\R^3$.

\begin{equation}\label{eq:NS}
\left\{\begin{aligned}
&\partial_t\rho^\eps+\Div(\rho^\eps u^\eps)=0,\\
&\partial_t(\rho^\eps u^\eps)
 +\Div(\rho^\eps u^\eps\otimes u^\eps)+\nabla p^\eps
 =\eps\Div\mathcal T(u^\eps),\\
&\partial_t(\rho^\eps E^\eps)
 +\Div\bigl((\rho^\eps E^\eps+p^\eps)u^\eps\bigr)
 =\eps\Div\bigl(\mathcal T(u^\eps)u^\eps
                  +\kappa\nabla\theta^\eps\bigr).
\end{aligned}\right.
\end{equation}
Here $\rho^\eps$, $u^\eps$ and $\theta^\eps$ denote the
density, velocity and temperature.
The pressure, specific total energy and viscous stress are
\begin{gather}
p=R\rho\theta,\qquad
E=\frac{R\theta}{\gamma-1}+\frac{|u|^2}{2},
\qquad R>0,\quad\gamma>1,\label{eq:gas}\\
\mathcal T(u)=\mu(\nabla u+\nabla u^{\mathsf T})
 +\lambda(\Div u)I,\qquad
\mu>0,\quad 2\mu+3\lambda\ge0,\quad\kappa>0.
\label{eq:coefficients}
\end{gather}
The transport coefficients are fixed, so viscosity and heat
conductivity vanish at the same rate.
We write
\[
v=\rho^{-1},\qquad q=(v,u,\theta),\qquad
U=(\rho,\rho u,\rho E)^{\mathsf T}.
\]
The last component, $U_E=\rho E$, is the total energy per unit volume.

The stability of viscous shocks has been studied extensively in
one space dimension.
For a one-dimensional model of a compressible viscous gas,
Matsumura--Nishihara \cite{MatsumuraNishihara1985} proved stability
under a zero-mass condition.
Goodman \cite{Goodman1986} developed a nonlinear stability theory
for viscous profiles of systems of conservation laws.
For the compressible Navier--Stokes equations, Liu
\cite{Liu1986} established stability by introducing a shift and
diffusion waves.
For the full system, Huang--Matsumura
\cite{HuangMatsumura2009} treated a composite wave consisting of
two small viscous shocks of comparable strengths.
Their result permits perturbations with nonzero mass.

Weighted relative entropy has also been used to study shock
stability. Kang--Vasseur \cite{KangVasseur2021a} proved a
contraction property for small barotropic viscous shocks under
large perturbations. They also established uniqueness and
stability of small isentropic Euler shocks in a class of inviscid
limits for density-dependent viscosities
\cite{KangVasseur2021b}. For the one-dimensional
Navier--Stokes--Fourier system, Kang--Vasseur--Wang
\cite{KangVasseurWang2023} introduced time-dependent shifts to
prove time-asymptotic stability for generic Riemann solutions
consisting of a viscous shock, a viscous contact wave and a
rarefaction wave.

For the one-dimensional vanishing-dissipation limit, Hoff--Liu
\cite{HoffLiu1989} treated shock data for the isentropic
Navier--Stokes equations. By combining matched asymptotic
expansions with energy estimates, Goodman--Xin
\cite{GoodmanXin1992} proved convergence for piecewise smooth
solutions with finitely many noninteracting shocks. For the
one-dimensional heat-conducting system, Wang \cite{Wang2008}
and Ma \cite{Ma2009} constructed viscous approximations for
prepared smooth data. Huang--Wang \cite{HuangWang2026}
obtained global existence for discontinuous initial data containing
a small shock. Their shifted relative-entropy estimate yields
convergence away from the shock on every finite time interval.

In several dimensions, Gu\`es--M\'etivier--Williams--Zumbrun
\cite{GuesMetivierWilliamsZumbrun2004} first constructed
viscous approximations to curved inviscid shocks for
parabolic regularizations. They later treated the physical
Navier--Stokes regularization of multidimensional Euler shocks
\cite{GuesMetivierWilliamsZumbrun2006}. These results hold
on the lifespan of the prescribed inviscid shock and rely
on structural and geometric hypotheses, together with an
Evans-function condition for spectral stability. For scalar
viscous conservation laws, Kang--Vasseur--Wang
\cite{KangVasseurWang2019} proved an $L^2$ contraction
estimate for large planar shocks with a space- and
time-dependent shift. For gas dynamics, Humpherys--Lyng--Zumbrun
\cite{HumpherysLyngZumbrun2017} studied the multidimensional
spectral stability of large-amplitude Navier--Stokes shocks
by asymptotic analysis and numerical Evans-function computations.

Wang--Wang \cite[Theorem~2.3]{WangWang2026} established the
nonlinear stability of weak planar shocks for the three-dimensional
isentropic Navier--Stokes equations on $\R\times\mathbb T^2$
under $H^2$ perturbations. Their proof combines weighted relative
entropy, higher-order energy estimates and a time-dependent shift,
without requiring a zero-mass condition. Chen--Cheng--Vasseur--Wang
\cite{ChenChengVasseurWang2026} proved that planar Euler shocks are
selected as vanishing-viscosity limits of weak solutions to the
two-dimensional barotropic Navier--Stokes equations on
$\R\times\mathbb T$ with density-dependent viscosity. Here we consider
the full heat-conducting system in $\R^3$ with constant transport
coefficients.

On the whole space, the estimates must also control the low frequencies.
Guo--Wang \cite{GuoWang2012} used negative Sobolev norms in
energy estimates for dissipative systems, including the compressible
Navier--Stokes equations. For planar shocks of the
full Navier--Stokes--Fourier system in $\R^3$, Huang--Xu
\cite[Theorem~1.1]{HuangXu2026} proved global nonlinear stability
and time-decay estimates for small conservative perturbations in
$H^3\cap\dot H^{-1}$. Their entropy-weighted elliptic energy,
constructed from the inverse entropy Hessian and a modified curl,
gives a time-integrated $L^2$ estimate. With this elliptic structure,
we close the estimates for $H^2\cap\dot H^{-1}$ initial data without
an $H^3$ bound on the perturbation.

We prove global stability of the viscous profile corresponding to a
weak planar Euler shock for the full compressible Navier--Stokes
equations in $\R^3$. For perturbations depending on all three spatial
variables, only smallness in a scaled $H^2\cap\dot H^{-1}$ norm is
required. The energy estimate holds uniformly for all $t\ge0$ and
gives the sharp $O(\sqrt\eps)$ convergence rate on every cylinder
$\R\times K$, where $K\subset\R^2$ has finite positive area.
If the scaled initial norm tends to zero, convergence away from the
shock front is uniform in time. The entropy-production measures
converge locally weakly to the measure carried by the Euler shock.

Three difficulties arise in closing the energy estimates.
\begin{enumerate}
\renewcommand{\labelenumi}{\arabic{enumi}.}
\item
The profile terms in the entropy identity have no definite sign,
and derivative estimates alone do not control the low frequencies
on $\R^3$. We therefore work in the stationary shock frame and use
an elliptic operator determined by the inverse entropy Hessian and
the profile. The corresponding negative-order energy gives a
time-integrated $L^2$ bound.

\item
At second order, the pressure term couples the highest density and
velocity derivatives. A direct estimate would lose one derivative.
Integration by parts moves the highest derivative from density to
velocity, where viscous dissipation controls it.

\item
Differentiating the viscous heating term twice produces second
derivatives of the reciprocal density multiplied by velocity
gradients. We integrate by parts before estimating these products,
leaving only first derivatives of the reciprocal density.
The third derivatives of velocity and temperature are then
controlled by the time-integrated second-order dissipation.
\end{enumerate}

The uniform $H^2$ estimate gives a common lifespan for the smooth
approximations obtained by Fourier truncation. Uniqueness follows
from a difference estimate, and the transport equation gives strong
continuity of the density perturbation in $H^2$.
The global energy estimate extends the local solution for all time.

After rescaling time and all three spatial variables, $A_{\eps,2}$
is exactly the $H^2\cap\dot H^{-1}$ norm of the rescaled conservative
initial perturbation. The unit-viscosity stability theorem therefore
applies to the rescaled problem. On each cylinder, the leading term
in the $L^2$ error is the exact distance between the viscous
profile and the Euler shock, yielding the sharp rate. Convergence
away from the shock front is uniform in time, including at $t=0$,
when $A_{\eps,2}\to0$. The same estimates yield local weak
convergence of the entropy-production measures to the shock measure.

Section~\ref{sec:main} gives the planar profile and the main
results. The proofs of stability, the vanishing-dissipation limit
and convergence of the entropy-production measures are given in
Section~\ref{sec:proof}. The uniform-in-time $H^2$ estimates
are established in Section~\ref{sec:H2}.

\section{Planar viscous shock and main results}
\label{sec:main}

\subsection{Planar viscous shock}

Let $q_\pm=(v_\pm,u_{1\pm}e_1,\theta_\pm)$ be the end
states of a small Lax 3-shock.
We assume $v_\pm,\theta_\pm>0$ and set
$e_1=(1,0,0)^{\mathsf T}$.
The shock speed is denoted by $\sigma$.
We write $\rho_\pm=v_\pm^{-1}$.
The jump of a quantity across the shock is denoted by
$[f]=f_+-f_-$.
The Rankine--Hugoniot relations are
\begin{equation}\label{eq:RH}
\left\{\begin{aligned}
&-\sigma[\rho]+[\rho u_1]=0,\\
&-\sigma[\rho u_1]+[\rho u_1^2+p]=0,\\
&-\sigma[\rho E]+[(\rho E+p)u_1]=0.
\end{aligned}\right.
\end{equation}
The Lax inequalities are
\[
 u_{1+}+\sqrt{\gamma R\theta_+}<\sigma
                 <u_{1-}+\sqrt{\gamma R\theta_-}.
\]
We denote the shock strength by $\delta=\rho_--\rho_+>0$
and set $\nu=2\mu+\lambda>0$.
The mass jump condition and the Lax inequality at the right
state imply
\[
 \rho_-(u_{1-}-\sigma)=\rho_+(u_{1+}-\sigma)<0.
\]
Substituting
$u_{1\pm}-\sigma=\rho_-(u_{1-}-\sigma)v_\pm$
into the momentum jump condition gives
\[
 [u_1]=\rho_-(u_{1-}-\sigma)[v],\qquad
 [p]=-\rho_-^2(u_{1-}-\sigma)^2[v].
\]
Substituting these relations into the energy jump condition, we obtain
\[
\begin{aligned}
 0&=\rho_-(u_{1-}-\sigma)[E]+[pu_1]\\
  &=\frac{R\rho_-(u_{1-}-\sigma)}{\gamma-1}[\theta]
       +\frac{u_{1+}+u_{1-}}2
          \left\{\rho_-(u_{1-}-\sigma)[u_1]+[p]\right\}\\
  &\quad+\frac{p_++p_-}{2}[u_1]\\
  &=\rho_-(u_{1-}-\sigma)
       \left\{\frac{R}{\gamma-1}[\theta]
                         +\frac{p_++p_-}{2}[v]\right\}.
\end{aligned}
\]
The momentum jump condition cancels the term containing
$u_{1+}+u_{1-}$. The density jump and negative mass flux give
$[v]=\delta/(\rho_-\rho_+)>0$
and $[u_1]<0$.

For $\eps=1$, introduce the shock coordinate $z=x_1-\sigma t$.
In this coordinate, write the viscous profile as
\[
 \widetilde q(z)
       =(\widetilde v(z),\widetilde u_1(z)e_1,\bar\theta(z)),
 \qquad \widetilde\rho=\widetilde v^{-1},\qquad
 \widetilde E=\frac{R\bar\theta}{\gamma-1}
                         +\frac12\widetilde u_1^2.
\]
A prime denotes differentiation with respect to $z$.
The profile equations in this coordinate are
\begin{equation}\label{eq:profile-ODE}
\left\{\begin{aligned}
&-\sigma\widetilde\rho'
                    +(\widetilde\rho\widetilde u_1)'=0,\\
&-\sigma(\widetilde\rho\widetilde u_1)'
   +(\widetilde\rho\widetilde u_1^2
                     +R\widetilde\rho\bar\theta)'
                                    =\nu\widetilde u_1'',\\
&-\sigma(\widetilde\rho\widetilde E)'
   +\{(\widetilde\rho\widetilde E
                        +R\widetilde\rho\bar\theta)\widetilde u_1\}'
           =(\nu\widetilde u_1\widetilde u_1'+\kappa\bar\theta')',
\end{aligned}\right.
\end{equation}
with $\widetilde q(\pm\infty)=q_\pm$.
Integrating the profile equations from $-\infty$ to $z$, we obtain
\[
\begin{aligned}
\widetilde\rho(\widetilde u_1-\sigma)
 &=\rho_-(u_{1-}-\sigma)=\rho_+(u_{1+}-\sigma)<0,\\
\nu\widetilde u_1'
 &=\rho_-(u_{1-}-\sigma)(\widetilde u_1-u_{1-})
       +R(\widetilde\rho\bar\theta-\rho_-\theta_-),\\
\nu\widetilde u_1\widetilde u_1'+\kappa\bar\theta'
 &=\rho_-(u_{1-}-\sigma)(\widetilde E-E_-)
       +R\widetilde\rho\bar\theta\widetilde u_1
       -R\rho_-\theta_-u_{1-}.
\end{aligned}
\]
Multiplying the second identity by $\widetilde u_1$ and
subtracting it from the third cancels the viscous-work term
and gives
\[
\begin{aligned}
 \kappa\bar\theta'
 &=\rho_-(u_{1-}-\sigma)
       \left\{\frac{R}{\gamma-1}(\bar\theta-\theta_-)
        +\frac{\widetilde u_1^2-u_{1-}^2}{2}
        -\widetilde u_1(\widetilde u_1-u_{1-})\right\}\\
 &\quad+R\rho_-\theta_-(\widetilde u_1-u_{1-})\\
 &=\frac{R\rho_-(u_{1-}-\sigma)}{\gamma-1}
                           (\bar\theta-\theta_-)
       +R\rho_-\theta_-(\widetilde u_1-u_{1-})\\
 &\quad-\frac{\rho_-(u_{1-}-\sigma)}2
                           (\widetilde u_1-u_{1-})^2.
\end{aligned}
\]

The changes $\theta\mapsto R\theta$ and $\kappa\mapsto\kappa/R$,
with $(u_{1-}+u_{1+})/2$ subtracted from the velocity and shock
speed, reduce the profile system to
\cite[equation (2.37)]{HuangXu2026}. By
\cite[Lemma~2.1]{HuangXu2026}, there is a smooth small profile,
unique up to translation. We use the original physical variables
here; the perturbation equations are normalized in
Section~\ref{sec:reformulation}.
We fix the translation by requiring
\[
 \widetilde u_1(0)=\frac{u_{1-}+u_{1+}}2.
\]
Let $\widetilde U$ denote the corresponding conservative
profile.
The end states and the chosen normalization do not depend on
$\eps$.
The following lemma records the profile estimates used below.

\begin{lemma}\label{lm:profile}
For sufficiently small $\delta$, the profile satisfies
\begin{equation}\label{g2:profile-bounds}
\begin{gathered}
 \widetilde u_1'<0,\qquad
 |\widetilde u_1'(z)|\le C\delta^2e^{-c\delta|z|},\\
 |\widetilde\rho'|+|\bar\theta'|\le C|\widetilde u_1'|,\qquad
 |(\widetilde\rho,\widetilde u_1,\bar\theta)''|
                  \le C\delta|\widetilde u_1'|,\\
 |\widetilde q(z)-q_\pm|\le C\delta e^{-c\delta|z|}
                       \quad\text{for }\ \pm z\ge0.
\end{gathered}
\end{equation}
\end{lemma}

\begin{proof}
The derivative bounds for
$(\widetilde\rho,\widetilde u_1,\bar\theta)$ follow from
\cite[Lemma~2.1, equations (2.38)--(2.39)]{HuangXu2026},
after applying the normalization above.
The Lax inequality at the right state gives
$\rho_-(\sigma-u_{1-})=\rho_+(\sigma-u_{1+})
 >\rho_+\sqrt{\gamma R\theta_+}$,
so the positive-state bounds keep the mass flux bounded away
from zero. The first profile equation then gives
\[
 \widetilde v'
   =\frac{\widetilde u_1'}{\rho_-(u_{1-}-\sigma)}>0,
 \qquad
 \widetilde\rho'=-\widetilde\rho^2\widetilde v'<0.
\]
The upper and lower bounds for $\widetilde\rho$, the estimate for
$\bar\theta'$ and the preceding identity imply
$|\widetilde q'|\le C|\widetilde u_1'|$. Integrating this bound
from $z$ to the corresponding end state and using the decay of
$\widetilde u_1'$, we obtain
\[
 |\widetilde q(z)-q_\pm|
 \le C\delta^2\int_{|z|}^\infty e^{-c\delta s}\dd s
 \le C\delta e^{-c\delta|z|}.
\]
\end{proof}

The higher-derivative estimate
\cite[equation (2.40)]{HuangXu2026} also gives
\[
 \norm{\partial_z^k
          (\widetilde\rho,\widetilde u_1,\bar\theta)}_\infty
       \le C_k\delta^2,\qquad 1\le k\le4.
\]
The coefficients below are smooth functions of
$(\widetilde\rho,\widetilde u_1,\bar\theta)$ on a fixed
neighborhood of the profile range. For instance,
\[
 \widetilde v'
       =-\widetilde\rho^{-2}\widetilde\rho',\qquad
 \widetilde v''
       =2\widetilde\rho^{-3}(\widetilde\rho')^2
                           -\widetilde\rho^{-2}\widetilde\rho''.
\]
The profile bounds and the chain rule give a
$C_k\delta^2$ bound for every coefficient derivative through
order four. Since the coefficients depend only on $x_1$,
Leibniz' rule gives
\[
 \partial^\alpha(\widetilde v f)
 =\sum_{j=0}^{\alpha_1}\binom{\alpha_1}{j}
       \partial_1^j(\widetilde v)
                        \partial^{\alpha-je_1}f.
\]
Taking $L^2$ norms and summing over $|\alpha|\le k$, we obtain
$\norm{\widetilde v f}_{H^k}\le C_k\norm f_{H^k}$ from the
coefficient bounds. The same multiplier estimate holds for every profile coefficient
used below. The constants $c,C,C_k$ may depend on the gas, the
transport coefficients and the fixed positive state neighborhood,
but are independent of $\eps$ and uniform for sufficiently small
$\delta$.

For $\eps>0$, the viscous profile and its Euler-shock
limit are given by
\begin{equation}\label{eq:profiles}
\widetilde q^\eps(t,x)=
 \widetilde q\left(\frac{x_1-\sigma t}{\eps}\right),
\qquad
q^s(t,x)=\begin{cases}
q_-,&x_1<\sigma t,\\
q_+,&x_1>\sigma t.
\end{cases}
\end{equation}
The conservative profile $\widetilde U^\eps$ has the same scaling.
Writing $z=(x_1-\sigma t)/\eps$, we compute
\[
\begin{gathered}
 \partial_t\widetilde q^\eps
       =-\frac\sigma\eps\widetilde q'(z),\qquad
 \partial_1^k\widetilde q^\eps
       =\eps^{-k}\widetilde q^{(k)}(z),\quad 1\le k\le4,\\
 \partial_2\widetilde q^\eps=\partial_3\widetilde q^\eps=0,\qquad
 \mathcal T(\widetilde u^\eps)
       =\eps^{-1}\widetilde u_1'(z)
                         \operatorname{diag}(\nu,\lambda,\lambda).
\end{gathered}
\]
The viscous and heat fluxes are therefore
\[
 \mathcal T(\widetilde u^\eps)\widetilde u^\eps
                      +\kappa\nabla\bar\theta^\eps
 =\eps^{-1}\bigl(\nu\widetilde u_1\widetilde u_1'
                              +\kappa\bar\theta'\bigr)(z)e_1.
\]
With these identities, \eqref{eq:NS} reduces to
\eqref{eq:profile-ODE} multiplied by $\eps^{-1}$, so
$\widetilde q^\eps$ is an exact solution of \eqref{eq:NS}.

\subsection{Global stability at unit viscosity}

All spatial norms below are taken over $\R^3$ unless stated
otherwise. For a five-component conservative perturbation $f$,
define
\begin{equation}\label{eq:negative-norm}
\norm f_{\dot H^{-1}}
 =\sup_{\substack{\varphi\in C_c^\infty(\R^3;\R^5)\\
                         \norm{\nabla\varphi}_2=1}}
       \abs{\langle f,\varphi\rangle}.
\end{equation}
For $k\ge0$, $\norm{\nabla^k f}_2^2$ is the sum of squared
derivative norms over all multi-indices of length $k$. We set
$\norm f_{H^k}^2=\sum_{j=0}^k\norm{\nabla^j f}_2^2$.

Relative to the fixed profile, the conservative and primitive
perturbations for \eqref{eq:NS} with $\eps=1$ are
\begin{equation}\label{g2:unit-perturbations}
\begin{gathered}
f(t,x)=U(t,x)-\widetilde U(x_1-\sigma t),\qquad
f_0(x)=U_0(x)-\widetilde U(x_1),\\
\Phi(t,x)=\log\frac{\rho(t,x)}
                         {\widetilde\rho(x_1-\sigma t)},\\
\psi(t,x)=u(t,x)-\widetilde u_1(x_1-\sigma t)e_1,\qquad
\zeta(t,x)=\theta(t,x)-\bar\theta(x_1-\sigma t),\quad
Z=(\Phi,\psi,\zeta).
\end{gathered}
\end{equation}

\begin{theorem}[Global $H^2$ stability]\label{g2:stability}
Fix the gas law and the transport coefficients in
\eqref{eq:gas}--\eqref{eq:coefficients}. For a sufficiently weak
planar Lax 3-shock, there exist constants $a_{*,2}>0$ and
$C_\delta>0$, independent of time, such that if the conservative
initial perturbation satisfies
\[
 f_0\in H^2(\R^3)\cap\dot H^{-1}(\R^3),\qquad
 \norm{f_0}_{H^2}^2+\norm{f_0}_{\dot H^{-1}}^2\le a_{*,2}^2,
\]
then the Cauchy problem for \eqref{eq:NS} with $\eps=1$ and
initial data $U_0$ admits a unique global strong solution.
The density and temperature remain uniformly bounded above
and away from zero. The perturbations have the regularity
\[
\begin{gathered}
 Z,f\in C([0,\infty);H^2(\R^3)),\qquad
 f\in C([0,\infty);\dot H^{-1}(\R^3)),\\
 \nabla(\psi,\zeta)\in L^2(0,\infty;H^2(\R^3)),
\end{gathered}
\]
and satisfy the energy estimate
\begin{equation}\label{g2:stability-energy}
\begin{aligned}
&\sup_{t\ge0}\left\{
       \norm{f(t)}_{\dot H^{-1}}^2+\norm{Z(t)}_{H^2}^2\right\}\\
&\quad+\int_0^\infty\left\{
       \norm Z_{H^2}^2+\norm{\nabla(\psi,\zeta)}_{H^2}^2
                         \right\}\dd t\\
&\qquad\le C_\delta\left\{
       \norm{f_0}_{\dot H^{-1}}^2+\norm{f_0}_{H^2}^2\right\}.
\end{aligned}
\end{equation}
Moreover, the fixed planar viscous profile is asymptotically
stable in the sense that
\[
 \lim_{t\to\infty}\left(
       \norm{Z(t)}_{H^2}+\norm{f(t)}_{H^2}\right)=0.
\]
The constants may depend on the fixed gas and transport
coefficients, the positive state neighborhood and the shock.
\end{theorem}

\subsection{Vanishing dissipation limit}

For conservative initial data
$U_0^\eps$, set
\begin{equation}\label{g2:preparation}
\begin{aligned}
w_0^\eps(x)&=U_0^\eps(x)-\widetilde U(x_1/\eps),\\
A_{\eps,2}^2
 &=\eps^{-5}\norm{w_0^\eps}_{\dot H^{-1}}^2
   +\sum_{k=0}^2\eps^{2k-3}\norm{\nabla^k w_0^\eps}_2^2.
\end{aligned}
\end{equation}
The powers of $\eps$ in \eqref{g2:preparation} are determined
by the change of variables $x=\eps y$. For $0\le k\le2$,
we have
\[
\begin{aligned}
 \norm{\nabla_y^k\{w_0^\eps(\eps\,\cdot)\}}_2^2
 &=\eps^{2k}\int_{\R^3}
                       |\nabla_x^kw_0^\eps(\eps y)|^2\dd y\\
 &=\eps^{2k-3}\norm{\nabla_x^kw_0^\eps}_2^2.
\end{aligned}
\]
The negative norm scales in the same way. In fact,
$\norm{\nabla\{\varphi(\cdot/\eps)\}}_2
       =\eps^{1/2}\norm{\nabla\varphi}_2$.
Thus $\varphi\mapsto\eps^{-1/2}\varphi(\cdot/\eps)$
maps the unit-gradient test-function class bijectively onto
itself. The definition \eqref{eq:negative-norm} then gives
\[
\begin{aligned}
 \norm{w_0^\eps(\eps\,\cdot)}_{\dot H^{-1}}
 &=\sup_{\norm{\nabla\varphi}_2=1}
       \eps^{-3}\abs{\langle w_0^\eps,
                                      \varphi(\cdot/\eps)\rangle}\\
 &=\eps^{-3}\eps^{1/2}\norm{w_0^\eps}_{\dot H^{-1}}\\
 &=\eps^{-5/2}\norm{w_0^\eps}_{\dot H^{-1}}.
\end{aligned}
\]
It follows that
\[
 A_{\eps,2}^2
 =\norm{w_0^\eps(\eps\,\cdot)}_{H^2}^2
       +\norm{w_0^\eps(\eps\,\cdot)}_{\dot H^{-1}}^2.
\]
In particular, the negative-order norm is imposed directly on
the conservative perturbation $w_0^\eps$. We use primitive
variables only in the nonnegative Sobolev norms.

For any solution, set
\begin{equation}\label{eq:perturbations}
\begin{gathered}
w^\eps=U^\eps-\widetilde U^\eps,\qquad
\psi^\eps=u^\eps-\widetilde u^\eps,\qquad
\zeta^\eps=\theta^\eps-\bar\theta^\eps,\\
Z^\eps=\left(\log\frac{\rho^\eps}{\widetilde\rho^\eps},
             \psi^\eps,\zeta^\eps\right),
\qquad \widetilde\rho^\eps=(\widetilde v^\eps)^{-1}.
\end{gathered}
\end{equation}
The following identities relate the density perturbations:
\[
\begin{aligned}
 \log\frac{\rho^\eps}{\widetilde\rho^\eps}
 &=(\rho^\eps-\widetilde\rho^\eps)
       \int_0^1\frac{\dd s}
        {\widetilde\rho^\eps+s(\rho^\eps-\widetilde\rho^\eps)},\\
 v^\eps-\widetilde v^\eps
 &=-\frac{\rho^\eps-\widetilde\rho^\eps}
                              {\rho^\eps\widetilde\rho^\eps}.
\end{aligned}
\]
The density bounds and these identities imply the pointwise
equivalence of $q^\eps-\widetilde q^\eps$ and $Z^\eps$ as long
as both states remain in a fixed positive neighborhood:
\begin{equation}\label{eq:pointwise-equivalence}
 C^{-1}|q^\eps-\widetilde q^\eps|
   \le |Z^\eps|\le C|q^\eps-\widetilde q^\eps|.
\end{equation}
We regard $q^s(0,\cdot)$ as a function of $x_1$ and set
\begin{equation}\label{eq:layer-constant}
 M_\delta=\norm{\widetilde q-q^s(0,\cdot)}_{L^2(\R)}
       \in(0,\infty).
\end{equation}
The decay estimate \eqref{g2:profile-bounds} yields
\[
\begin{aligned}
 M_\delta^2
 &=\int_{-\infty}^0|\widetilde q(z)-q_-|^2\dd z
       +\int_0^\infty|\widetilde q(z)-q_+|^2\dd z\\
 &\le C\delta^2\int_\R e^{-2c\delta|z|}\dd z
 \le C\delta.
\end{aligned}
\]
Thus $M_\delta<\infty$. If $M_\delta=0$, then
$\widetilde q=q_-$ on $(-\infty,0)$ and
$\widetilde q=q_+$ on $(0,\infty)$. Continuity at $z=0$
would imply $q_-=q_+$, contrary to $\delta>0$.

\begin{corollary}[Vanishing dissipation limit]\label{g2:theorem}
Let the gas law, transport coefficients and weak planar shock
satisfy the assumptions of Theorem~\ref{g2:stability}, and let
$a_{*,2}$ be the smallness constant there. Then there exist constants
$c_\delta,C_\delta>0$, independent of $0<\eps\le1$, such that,
if $A_{\eps,2}\le a_{*,2}$, the Cauchy problem for
\eqref{eq:NS} with initial data $U_0^\eps$ admits a unique
global strong solution. The density and temperature remain uniformly bounded above
and away from zero for all time. For each fixed $\eps$, the
perturbations have the regularity
\[
\begin{gathered}
Z^\eps,w^\eps\in C([0,\infty);H^2(\R^3)),\qquad
w^\eps\in C([0,\infty);\dot H^{-1}(\R^3)),\\
\nabla(\psi^\eps,\zeta^\eps)
             \in L^2(0,\infty;H^2(\R^3)).
\end{gathered}
\]
The solution satisfies the uniform energy estimate
\begin{equation}\label{g2:physical-energy}
\begin{aligned}
&\sup_{t\ge0}\left\{
 \eps^{-2}\norm{w^\eps(t)}_{\dot H^{-1}}^2
   +\sum_{k=0}^2\eps^{2k}\norm{\nabla^k Z^\eps(t)}_2^2\right\}\\
&\quad+\int_0^\infty\left\{
 \sum_{k=0}^2\eps^{2k-1}\norm{\nabla^k Z^\eps}_2^2
 +\sum_{k=0}^2\eps^{2k+1}
           \norm{\nabla^{k+1}(\psi^\eps,\zeta^\eps)}_2^2
                       \right\}\dd t
 \le C_\delta\eps^3 A_{\eps,2}^2.
\end{aligned}
\end{equation}
For every measurable set $K\subset\R^2$ with
$0<|K|<\infty$, we have
\begin{equation}\label{g2:cylinder}
\sup_{t\ge0}
\left|\norm{q^\eps(t)-q^s(t)}_{L^2(\R\times K)}
       -|K|^{1/2}M_\delta\sqrt\eps\right|
 \le C_\delta\eps^{3/2}A_{\eps,2}.
\end{equation}
Thus the rate $O(\sqrt\eps)$ is sharp on every such
cylinder. For each $\ell>0$, the error away from the shock
front satisfies
\begin{equation}\label{g2:away}
\sup_{t\ge0}
\norm{q^\eps(t)-q^s(t)}_{L^\infty(\{|x_1-\sigma t|\ge\ell\})}
 \le C_\delta A_{\eps,2}+C_\delta e^{-c_\delta\ell/\eps}.
\end{equation}
When $A_{\eps,2}\to0$, estimate \eqref{g2:away} gives
convergence uniformly in time outside any fixed neighborhood
of the shock front, including at $t=0$. For any family satisfying
$A_{\eps,2}\le a_{*,2}$, the weak formulations of
\eqref{eq:NS} converge to the Euler equations with state
$q^s$. The corresponding entropy production measures
converge locally weakly to the measure in
\eqref{eq:production-limit}.
\end{corollary}

Sobolev embedding in the rescaled variables gives
\[
 \norm{w_0^\eps}_\infty
 =\norm{w_0^\eps(\eps\,\cdot)}_\infty
 \le C\left(\sum_{k=0}^2
          \eps^{2k-3}\norm{\nabla^k w_0^\eps}_2^2\right)^{1/2}
 \le C A_{\eps,2}.
\]
The profile takes values in a compact subset of the positive
state space. For sufficiently small $A_{\eps,2}$, the preceding
bound therefore keeps the initial density and temperature positive.
Perturbations may depend on all three spatial variables, with
neither an $L^1$ bound nor a zero-mass condition. The constants
may depend on the gas and transport coefficients, the fixed
positive state neighborhood and the shock; they are independent
of $\eps$ and time.

\begin{remark}\label{rem:global-L2}
The equations and the energy estimate are posed on $\R^3$.
The convergence estimate \eqref{g2:cylinder} holds on each
cylinder with a transverse section of finite positive area.
Since the profile difference $\widetilde q^\eps-q^s$ does not
depend on $x'$, Tonelli's theorem and the substitution
$x_1=\sigma t+\eps z$ give
\[
\begin{aligned}
 \norm{\widetilde q^\eps(t)-q^s(t)}_{L^2(\R^3)}^2
 &=\int_{\R^2}\eps\int_\R
             |\widetilde q(z)-q^s(0,z)|^2\dd z\dd x'\\
 &=\eps M_\delta^2\int_{\R^2}\dd x'=\infty.
\end{aligned}
\]
By \eqref{eq:pointwise-equivalence} and
\eqref{g2:physical-energy},
$q^\eps-\widetilde q^\eps\in L^2(\R^3)$.
The norm calculation above then shows that
$q^\eps-q^s$ cannot belong to $L^2(\R^3)$.
\end{remark}

\begin{corollary}\label{cor:example}
Let $G\in H^2(\R^3;\R^5)\cap\dot H^{-1}(\R^3;\R^5)$
and let $\beta>5/2$. For each viscosity parameter, choose
\begin{equation}\label{eq:example}
 U_0^\eps(x)=\widetilde U(x_1/\eps)+\eps^\beta G(x).
\end{equation}
For all sufficiently small $\eps$, the data in
\eqref{eq:example} satisfy the assumptions of
Corollary~\ref{g2:theorem}, and the corresponding solutions
converge uniformly away from the shock front, including at
the initial time.
The function $G$ need not be planar and may depend on all
three spatial variables.
\end{corollary}

\begin{proof}
Substituting \eqref{eq:example} into \eqref{g2:preparation}
gives
\[
 A_{\eps,2}^2
 =\eps^{2\beta-5}\norm G_{\dot H^{-1}}^2
  +\sum_{k=0}^2\eps^{2\beta+2k-3}\norm{\nabla^kG}_2^2.
\]
For $0<\eps\le1$ and $0\le k\le2$, we have
$\eps^{2\beta+2k-3}=\eps^{2\beta-5}\eps^{2k+2}
       \le\eps^{2\beta-5}$. Therefore
\[
 A_{\eps,2}^2
 \le\eps^{2\beta-5}
       \left(\norm G_{\dot H^{-1}}^2
                         +\sum_{k=0}^2\norm{\nabla^kG}_2^2\right)
 \longrightarrow0.
\]
Consequently, $A_{\eps,2}\le C_G\eps^{\beta-5/2}$,
so the required smallness condition holds for sufficiently
small viscosity.
Since the profile takes values in a compact subset of the
positive state space and $\eps^\beta\norm G_\infty\to0$,
the initial density and temperature are positive for
sufficiently small $\eps$.
The error terms containing $A_{\eps,2}$ in
\eqref{g2:cylinder} and \eqref{g2:away} satisfy
\[
 \eps^{3/2}A_{\eps,2}\le C_G\eps^{\beta-1},\qquad
 \eps^{-1/2}\eps^{3/2}A_{\eps,2}
                       \le C_G\eps^{\beta-3/2},\qquad
 A_{\eps,2}\le C_G\eps^{\beta-5/2}.
\]
Dividing the middle bound by $|K|^{1/2}M_\delta>0$
controls the error relative to the leading term in
\eqref{g2:cylinder}. The last bound tends to zero and
gives the asserted uniform convergence away from the shock
front, including at the initial time.
\end{proof}

We give an explicit example of such a function:
\[
 G(x)=(1+|x|^2)^{-4/3}(1,0,0,0,0)^{\mathsf T}
       \in H^2\cap L^{6/5}\subset H^2\cap\dot H^{-1},
 \qquad G\notin L^1.
\]
By \cite[Lemma~A.6]{GuoWang2012}, the Riesz potential of order
one maps $L^{6/5}(\R^3)$ to $L^2(\R^3)$. Thus
\[
 \norm G_{\dot H^{-1}}\le C\norm G_{6/5}.
\]
A direct calculation shows that
\[
 \left|\partial^\alpha(1+|x|^2)^{-4/3}\right|
 \le C_\alpha(1+|x|^2)^{-4/3-|\alpha|/2},
 \qquad |\alpha|\le2,
\]
Hence all derivatives of order at most two belong to
$L^2(\R^3)$. Indeed, in polar coordinates,
\[
\begin{aligned}
 \int_{\R^3}|G|^{6/5}\dd x
   &=4\pi\int_0^\infty r^2(1+r^2)^{-8/5}\dd r<\infty,\\
 \int_{\R^3}|G|\dd x
   &=4\pi\int_0^\infty r^2(1+r^2)^{-4/3}\dd r=\infty.
\end{aligned}
\]
The same duality estimate applies to every smooth compactly
supported function $G$. Such functions belong to
$H^2\cap L^{6/5}$ and may have nonzero spatial integral.

\begin{remark}\label{rem:H2-preparation}
The $H^2$ preparation permits smooth initial perturbations whose
scaled $H^3$ norm is arbitrarily large. Fix a nonzero real function
$\chi\in C_c^\infty(\R^3)$ and let
\[
 G_n(y)=a n^{-2}\chi(y)\sin(ny_1)
                       (1,0,0,0,0)^{\mathsf T},\qquad n\ge1.
\]
Leibniz's rule and the embedding
$L^{6/5}\hookrightarrow\dot H^{-1}$ give
\[
 \norm{G_n}_{H^2}\le Ca,\qquad
 \norm{G_n}_{\dot H^{-1}}
       \le C\norm{G_n}_{6/5}\le Ca n^{-2}.
\]
Direct differentiation gives
\[
 \partial_1^3G_n
 =-an\chi(y)\cos(ny_1)(1,0,0,0,0)^{\mathsf T}
       +O_{L^2}(a).
\]
Integration by parts gives
\[
 \norm{\chi\cos(ny_1)}_2^2
 =\frac12\norm\chi_2^2
       +\frac12\int\chi^2\cos(2ny_1)\dd y
 \longrightarrow\frac12\norm\chi_2^2>0,
\]
and hence $\norm{G_n}_{H^3}\sim an$ as $n\to\infty$.
The constants depend on $\chi$ and are independent of $a$ and $n$.
Choose $a>0$ sufficiently small and integers $n(\eps)\to\infty$,
and set
\[
 U_0^\eps(x)=\widetilde U(x_1/\eps)+G_{n(\eps)}(x/\eps).
\]
The scaling identity for the preparation gives
\[
 A_{\eps,2}^2
 =\norm{G_{n(\eps)}}_{H^2}^2
       +\norm{G_{n(\eps)}}_{\dot H^{-1}}^2\le Ca^2,
 \qquad
 \norm{w_0^\eps(\eps\,\cdot)}_{H^3}\sim a n(\eps)
       \longrightarrow\infty.
\]
Thus Corollary~\ref{g2:theorem} applies uniformly to this family.
Taking $a=a(\eps)\to0$ with $a(\eps)n(\eps)\to\infty$
also gives $A_{\eps,2}\to0$ and hence uniform convergence
away from the shock front, including at the initial time.
\end{remark}

\section{Proofs of the main results}
\label{sec:proof}

\subsection{Reformulation of the problem}
\label{sec:reformulation}

We consider \eqref{eq:NS} with $\eps=1$ and normalize the
temperature by
\[
 \theta\mapsto R\theta,\qquad c_v=\frac1{\gamma-1},\qquad
 \kappa\mapsto\frac\kappa R.
\]
We keep the notation $\theta$ for the rescaled temperature.
The pressure law and specific total energy become
\[
 p=\rho\theta,\qquad E=c_v\theta+\frac{|u|^2}{2}.
\]
This normalization leaves the conservative variable $U$
unchanged. We next apply the Galilean transformation
\[
 x\mapsto x-ct e_1,\qquad u\mapsto u-ce_1,\qquad
 \sigma\mapsto\sigma-c,\qquad c=\frac{u_{1-}+u_{1+}}2.
\]
The velocities at the two end states then have equal magnitude
and opposite signs. For the momentum and total energy densities,
the transformation gives
\[
\begin{aligned}
 \rho(u-ce_1)&=\rho u-c\rho e_1,\\
 \rho\left(c_v\theta+\frac{|u-ce_1|^2}{2}\right)
       &=\rho E-c\rho u_1+\frac{c^2}{2}\rho.
\end{aligned}
\]
A conservative perturbation $(r,f_m,f_E)^{\mathsf T}$ therefore
transforms according to
\begin{equation}\label{eq:galilean}
(r,f_m,f_E)\longmapsto
 \left(r,f_m-cr e_1,f_E-cf_{m,1}+\tfrac12c^2r\right).
\end{equation}
The transformations with parameters $c$ and $-c$ are inverse
to each other and bounded on $H^2$ and $\dot H^{-1}$.
Spatial translations preserve these norms. The two frames therefore
give equivalent norms, with constants independent of $\eps$. We keep
the notation $\kappa$ for the normalized heat conductivity and
$\sigma$ for the shock speed in the new frame.

The shock coordinates are
\[
 \tau=t,\qquad y=x-\sigma te_1
.\]
The profile is stationary, with
$\partial_t=\partial_\tau-\sigma\partial_{y_1}$ and
$\nabla_x=\nabla_y$. All spatial derivatives below are taken
with respect to $y$.

Let $\rho_*=(\rho_-+\rho_+)/2$ and
$\theta_*=(\theta_-+\theta_+)/2$. The profile stays within
$O(\delta)$ of $(\rho_*,0,\theta_*)$.
Write $\widetilde u=\widetilde u_1e_1$ and
$f=U-\widetilde U=(r,f_m,f_E)^{\mathsf T}$, with
$r=\rho-\widetilde\rho$. Define
\[
\begin{gathered}
 Z=(\Phi,\psi,\zeta),\qquad
 \Phi=\log(\rho/\widetilde\rho),\qquad
 \psi=u-\widetilde u,\qquad \zeta=\theta-\bar\theta,\\
 N=\norm Z_{H^2}.
\end{gathered}
\]
Throughout this subsection, we consider perturbations in a
small $H^2$ ball and assume that the density and temperature
remain positive. Expanding the conservative perturbation in the primitive
variables gives
\[
\begin{aligned}
 r&=\widetilde\rho(e^\Phi-1),\qquad
 f_m=\rho\psi+r\widetilde u,\\
 f_E&=c_v(\rho\zeta+\bar\theta r)
       +\frac12r|\widetilde u|^2
       +\rho\widetilde u\cdot\psi+\frac12\rho|\psi|^2.
\end{aligned}
\]
Solving these identities for the primitive perturbations gives
\[
\begin{aligned}
 \Phi&=\log(1+r/\widetilde\rho),\qquad
 \psi=\rho^{-1}(f_m-r\widetilde u),\\
 c_v\rho\zeta
   &=f_E-\left(c_v\bar\theta+\frac12|\widetilde u|^2\right)r
       -\rho\widetilde u\cdot\psi-\frac12\rho|\psi|^2.
\end{aligned}
\]
In particular,
$v-\widetilde v=-r/(\rho\widetilde\rho)$ vanishes with the
density perturbation. Differentiating the density relation gives
\[
\begin{aligned}
 \nabla r&=\rho\nabla\Phi
                 +\widetilde\rho'(e^\Phi-1)e_1,\\
 \nabla^2r&=\rho\nabla^2\Phi
                   +\rho\nabla\Phi\otimes\nabla\Phi\\
 &\quad+\widetilde\rho'e^\Phi
       (e_1\otimes\nabla\Phi+\nabla\Phi\otimes e_1)
       +\widetilde\rho''(e^\Phi-1)e_1\otimes e_1.
\end{aligned}
\]
Sobolev embedding and interpolation give
\[
\begin{aligned}
 \norm{\rho\nabla\Phi\otimes\nabla\Phi}_2
 &\le C\norm{\nabla\Phi}_3\norm{\nabla\Phi}_6\\
 &\le C\norm\Phi_{H^2}^2
 \le CN\norm\Phi_{H^2}.
\end{aligned}
\]
The inverse density relation gives
$\nabla\Phi=v\nabla r+(v-\widetilde v)\widetilde\rho'e_1$.
Applying Sobolev product estimates in both directions to
$r,f_m,f_E$, with the profile coefficients controlled by
Lemma~\ref{lm:profile}, we obtain
\begin{equation}\label{g2:positive-equivalence}
 \norm f_{H^j}\sim\norm Z_{H^j},\qquad 0\le j\le2.
\end{equation}
Set
$\mathcal L=\mu\Delta+(\mu+\lambda)\nabla\Div$.
With the normalizations above, the perturbation equations
\cite[equations (3.22)--(3.23)]{HuangXu2026} take the form
\begin{equation}\label{ft:primitive}
\left\{\begin{aligned}
&D_\tau\Phi+\Div\psi=F_\Phi,\\
&D_\tau\psi+\theta\nabla\Phi+\nabla\zeta-v\mathcal L\psi=F_\psi,\\
&c_vD_\tau\zeta+\theta\Div\psi-\kappa v\Delta\zeta=F_\zeta,
\end{aligned}\right.
\qquad D_\tau=\partial_\tau+(u-\sigma e_1)\cdot\nabla,
\end{equation}
where
\begin{align*}
F_\Phi&=-\frac{\widetilde\rho'}{\widetilde\rho}\psi_1,\\
F_\psi&=-\widetilde u_1'\psi_1e_1
 -\frac{\widetilde\rho'}{\widetilde\rho}\zeta e_1
 +\nu\widetilde u_1''(v-\widetilde v) e_1,\\
F_\zeta&=-c_v\bar\theta'\psi_1-\widetilde u_1'\zeta
 +\kappa(v-\widetilde v)\bar\theta''
 +\nu(\widetilde u_1')^2(v-\widetilde v)
 +v\mathcal T(\psi):\nabla\psi
 +2v\mathcal T(\widetilde u):\nabla\psi.
\end{align*}
Here $\Phi$, $v$ and $\mathcal T$ correspond to $\phi$, $a$ and
$S$ in the cited equations. All source terms vanish at $Z=0$.

For the second-order estimate, define the energy and dissipation
functionals
\begin{align}
S_2(Z)&=\frac12\sum_{|\alpha|\le2}\int_{\R^3}
 \left(\theta_*|\partial^\alpha\Phi|^2
       +|\partial^\alpha\psi|^2
       +\frac{c_v}{\theta_*}|\partial^\alpha\zeta|^2\right)\dd y,
 \nonumber\\
I_2(Z)&=\sum_{|\alpha|\le1}\int_{\R^3}
                  \partial^\alpha\psi\cdot\nabla\partial^\alpha\Phi\dd y,
 \qquad P_2=S_2+\sigma_2 I_2,\label{ft:P2}\\
D_2(Z)&=\norm{\nabla\Phi}_{H^1}^2
                    +\norm{\nabla(\psi,\zeta)}_{H^2}^2.\nonumber
\end{align}
Cauchy--Schwarz gives
\[
 |I_2(Z)|\le C\norm\psi_{H^1}\norm{\nabla\Phi}_{H^1}
          \le C\norm Z_{H^2}^2.
\]
Since $\theta_*,c_v>0$, the preceding bound allows us to
choose $\sigma_2>0$ so that $\sigma_2|I_2|\le S_2/2$.
Thus
\[
 \frac12S_2\le P_2\le\frac32S_2,
 \qquad P_2\sim\norm Z_{H^2}^2.
\]
Lemma~\ref{ft:H2-lemma} gives the additional restriction
on $\sigma_2$ required by the energy estimate.

\subsection{Local existence}

\begin{lemma}\label{g2:local}
Let $Z_0\in H^2(\R^3)$ be sufficiently small.
Then there exists $T>0$ such that the stationary perturbation
system \eqref{ft:primitive} with $Z(0)=Z_0$ admits a unique
solution with the regularity
\[
Z\in C([0,T];H^2),\qquad
\nabla(\psi,\zeta)\in L^2(0,T;H^2),\qquad
(\psi_\tau,\zeta_\tau)\in L^2(0,T;H^1).
\]
The density and temperature remain positive. For a sufficiently
small $b>0$, the solutions with $\norm{Z_0}_{H^2}\le b$
exist on a common positive time interval and remain in a
fixed larger small ball.
A solution can therefore be continued past any finite time
provided that $\sup_\tau\norm{Z(\tau)}_{H^2}\le b$.
If $f_0\in\dot H^{-1}$, then both $f$ and the effective
perturbation $g$ from \eqref{ft:effective} belong to
$H^1(0,T;\dot H^{-1})$.
\end{lemma}

\begin{proof}
For a related problem with constant far-field density,
Cho--Kim \cite[Section~3, Theorem~3.6 of the 2004 preprint]{ChoKim2006}
established local existence with $H^2$ regularity of the density.
Divide the last equation in \eqref{ft:primitive} by $c_v$
and write the system as $Z_\tau=\mathcal F(Z)$. Let $J_n$
be the orthogonal Fourier projection onto $\{|\xi|\le n\}$.
We solve the truncated system
\[
\partial_\tau Z_n=J_n\mathcal F(Z_n),\qquad Z_n(0)=J_nZ_0.
\]
The projection has the bounds
\[
 \norm{\nabla^kJ_nh}_{H^2}\le n^k\norm h_{H^2},
 \qquad
 \norm{J_nh}_{H^2}\le\norm h_{H^2}.
\]
Since $J_n$ projects onto frequencies in $\{|\xi|\le n\}$,
spatial differentiation is bounded on $\operatorname{Ran}J_n$.
The truncated equation is thus an ordinary differential
equation on this range. The multiplier bounds for the profile
coefficients show that its vector field is locally Lipschitz
in $H^2$. The profile equations also imply
$\mathcal F(0)=0$.

In the $S_2$ pairing, self-adjointness of $J_n$,
commutation with spatial derivatives and \mbox{$J_nZ_n=Z_n$} give
\[
 \left\langle J_n\partial^\alpha\mathcal F(Z_n),
                         \partial^\alpha Z_n\right\rangle
 =\left\langle\partial^\alpha\mathcal F(Z_n),
                         \partial^\alpha Z_n\right\rangle,
 \qquad |\alpha|\le2.
\]
The same identity holds for the cross term because
$\nabla\partial^\alpha\Phi_n$ and $\partial^\alpha\psi_n$
belong to $\operatorname{Ran}J_n$. Therefore,
\[
\begin{aligned}
 \frac{\dd}{\dd\tau}I_2(Z_n)
 &=\sum_{|\alpha|\le1}\bigl\{
 \langle\partial^\alpha\mathcal F_\psi(Z_n),
                         \nabla\partial^\alpha\Phi_n\rangle\\
 &\hspace{38mm}+\langle\partial^\alpha\psi_n,
                     \nabla\partial^\alpha\mathcal F_\Phi(Z_n)\rangle
                                                        \bigr\}.
\end{aligned}
\]
Lemma~\ref{ft:H2-lemma} gives an estimate uniform in $n$.
For sufficiently small $\delta+\norm{Z_n}_{H^2}$,
absorb the $D_2$ term in \eqref{ft:H2-energy} to obtain
\[
 \frac{\dd}{\dd\tau}P_2(Z_n)+cD_2(Z_n)
                       \le C_\delta P_2(Z_n).
\]
The equivalence $P_2(Z_n)\sim\norm{Z_n}_{H^2}^2$
and the contraction property
$\norm{J_nZ_0}_{H^2}\le\norm{Z_0}_{H^2}$ allow us
to apply Gronwall's inequality. We obtain
\[
 P_2(Z_n(\tau))+c\int_0^\tau D_2(Z_n(s))\dd s
       \le C_\delta e^{C_\delta\tau}\norm{Z_0}_{H^2}^2.
\]
Increase $C_\delta\ge1$ to include the constants in
$P_2(Z_n)\sim\norm{Z_n}_{H^2}^2$. Choose $b>0$ so that
the ball of radius $2C_\delta^{1/2}b$ is contained in the
small ball on which the estimate holds. We then choose
$T>0$ such that $e^{C_\delta T}\le2$.

As long as $\norm{Z_n}_{H^2}\le2C_\delta^{1/2}b$,
the initial bound $\norm{Z_0}_{H^2}\le b$ and the
preceding estimate imply
\[
 \norm{Z_n(\tau)}_{H^2}^2
 \le C_\delta e^{C_\delta\tau}\norm{Z_0}_{H^2}^2
 \le2C_\delta b^2<4C_\delta b^2,\qquad \tau\le T.
\]
This strictly improves the assumed bound. The ODE
continuation criterion therefore gives a common existence
interval $[0,T]$, on which
\begin{equation}\label{g2:local-bound}
\sup_{0\le\tau\le T}\norm{Z_n(\tau)}_{H^2}^2
       +\int_0^T D_2(Z_n)\dd\tau
 \le C_\delta\norm{Z_0}_{H^2}^2.
\end{equation}
On this ball, \eqref{ft:embeddings} and
\eqref{ft:heating} imply
\[
 \norm{\mathcal F(Z_n)}_{H^1}
 \le C_\delta\left(\norm{Z_n}_{H^2}
                 +\norm{\nabla(\psi_n,\zeta_n)}_{H^2}\right).
\]
For the diffusion terms, use
$v_n=\widetilde v+(v_n-\widetilde v)$, where
$v_n=\widetilde v e^{-\Phi_n}$. This gives
\[
\begin{aligned}
 &\norm{v_n\mathcal L\psi_n}_{H^1}
       +\norm{\kappa v_n\Delta\zeta_n}_{H^1}\\
 &\quad\le C_\delta\bigl(1+\norm{Z_n}_{H^2}\bigr)
                         \norm{\nabla(\psi_n,\zeta_n)}_{H^2}.
\end{aligned}
\]
For the viscous heating term, differentiation of $v_n$
produces a profile term and the perturbative term
$-v_n\nabla\Phi_n$. The product estimates in
\eqref{ft:heating} give
\[
\begin{aligned}
 &\norm{v_n\mathcal T(\psi_n):\nabla\psi_n}_{H^1}\\
 &\quad\le C_\delta\norm{\nabla\psi_n}_3\norm{\nabla\psi_n}_6
       +C\norm{\nabla\Phi_n}_6\norm{\nabla\psi_n}_6^2\\
 &\qquad\quad+C\norm{\nabla\psi_n}_6\norm{\nabla^2\psi_n}_3
 \le C_\delta\norm{Z_n}_{H^2}\norm{\nabla\psi_n}_{H^2}.
\end{aligned}
\]
The bound \eqref{g2:local-bound} now gives
\[
\norm{\mathcal F(Z_n)}_{L^2(0,T;H^1)}\le C_{\delta,T}.
\]
It follows that the truncation error
$R_n=(J_n-I)\mathcal F(Z_n)$ tends to zero
in $L^2(0,T;L^2)$; in fact,
\[
\norm{R_n}_{L^2(0,T;L^2)}
 \le n^{-1}\norm{\mathcal F(Z_n)}_{L^2(0,T;H^1)}.
\]
Thus the approximate equations take the form
$\partial_\tau Z_n=\mathcal F(Z_n)+R_n$.
To prove convergence, set
$\delta Z=Z_n-Z_m$, $\delta R=R_n-R_m$, and
$\delta v=v_n-v_m$.
The identity $v_n=\widetilde v e^{-\Phi_n}$
and smoothness of the reciprocal density on the state range give
$\norm{\delta v}_2\le C\norm{\delta\Phi}_2$.
We use the difference energy
\[
 S_0(\delta Z)=\frac12\int
 \left(\theta_*|\delta\Phi|^2+|\delta\psi|^2
                  +\frac{c_v}{\theta_*}|\delta\zeta|^2\right)\dd y.
\]
Writing $u_n=\widetilde u+\psi_n$ and
$\theta_n=\bar\theta+\zeta_n$, we find that $\delta Z=Z_n-Z_m$
satisfies
\[
\begin{aligned}
 &\partial_\tau\delta\Phi
       +(u_n-\sigma e_1)\cdot\nabla\delta\Phi+\Div\delta\psi\\
 &\qquad=-\delta\psi\cdot\nabla\Phi_m
       +F_\Phi(Z_n)-F_\Phi(Z_m)+(\delta R)_\Phi,\\
 &\partial_\tau\delta\psi
       +(u_n-\sigma e_1)\cdot\nabla\delta\psi
       +\theta_n\nabla\delta\Phi+\nabla\delta\zeta
       -v_n\mathcal L\delta\psi\\
 &\qquad=-\delta\psi\cdot\nabla\psi_m
       -\delta\zeta\nabla\Phi_m+\delta v\mathcal L\psi_m
       +F_\psi(Z_n)-F_\psi(Z_m)+(\delta R)_\psi,\\
 &c_v\partial_\tau\delta\zeta
       +c_v(u_n-\sigma e_1)\cdot\nabla\delta\zeta
       +\theta_n\Div\delta\psi-\kappa v_n\Delta\delta\zeta\\
 &\qquad=-c_v\delta\psi\cdot\nabla\zeta_m
       -\delta\zeta\Div\psi_m+\kappa\delta v\Delta\zeta_m\\
 &\qquad\quad+F_\zeta(Z_n)-F_\zeta(Z_m)+c_v(\delta R)_\zeta.
\end{aligned}
\]
Testing the equations with
$\theta_*\delta\Phi$, $\delta\psi$, and
$\delta\zeta/\theta_*$, respectively, and integrating over $\R^3$,
we find that the constant-coefficient acoustic terms cancel:
\[
\begin{aligned}
 \int\bigl(&\theta_*\delta\Phi\Div\delta\psi
       +\theta_*\nabla\delta\Phi\cdot\delta\psi\\
       &+\nabla\delta\zeta\cdot\delta\psi
                         +\delta\zeta\Div\delta\psi\bigr)\dd y=0.
\end{aligned}
\]
Since $S_0(\delta Z)\sim\norm{\delta Z}_2^2$,
integration by parts in the transport terms gives
\[
\begin{aligned}
 &\frac12\left|\int\Div u_n
       \left(\theta_*|\delta\Phi|^2+|\delta\psi|^2
                  +\frac{c_v}{\theta_*}|\delta\zeta|^2\right)\dd y\right|\\
 &\quad\le C\bigl(\delta^2+\norm{\nabla\psi_n}_\infty\bigr)
                                             \norm{\delta Z}_2^2\\
 &\quad\le C_\delta\bigl(1+D_2(Z_n)\bigr)\norm{\delta Z}_2^2.
\end{aligned}
\]
The variable diffusion coefficient satisfies
\[
 \nabla v_n=\widetilde v'e^{-\Phi_n}e_1
                                      -v_n\nabla\Phi_n.
\]
For the velocity diffusion term, integration by parts gives
\[
\begin{aligned}
 -\int v_n\mathcal L\delta\psi\cdot\delta\psi\dd y
 &=\mu\int v_n|\nabla\delta\psi|^2\dd y
       +(\mu+\lambda)\int v_n|\Div\delta\psi|^2\dd y\\
 &\quad+\mu\sum_{j=1}^3\int(\partial_jv_n)
                     \partial_j\delta\psi\cdot\delta\psi\dd y\\
 &\quad+(\mu+\lambda)\int(\nabla v_n\cdot\delta\psi)
                                         \Div\delta\psi\dd y.
\end{aligned}
\]
The coefficient $v_n$ has a positive lower bound independent
of $n$, and $\mu+\lambda\ge\mu/3>0$.
The first two terms are coercive and control $\norm{\nabla\delta\psi}_2^2$.
The remaining terms satisfy
\[
\begin{aligned}
 &C\int|\nabla v_n|\,|\nabla\delta\psi|\,|\delta\psi|\dd y\\
 &\quad\le C\delta^2\norm{\nabla\delta\psi}_2\norm{\delta\psi}_2
       +C\norm{\nabla\Phi_n}_3
                  \norm{\nabla\delta\psi}_2\norm{\delta\psi}_6\\
 &\quad\le\bigl(\eta_0+C\norm{Z_n}_{H^2}\bigr)
                       \norm{\nabla\delta\psi}_2^2
       +C_{\eta_0}\delta^4\norm{\delta\psi}_2^2.
\end{aligned}
\]
The temperature diffusion term is treated in the same way:
\[
 -\int v_n\Delta\delta\zeta\,\delta\zeta\dd y
 =\int v_n|\nabla\delta\zeta|^2\dd y
       +\int\delta\zeta\,\nabla v_n\cdot\nabla\delta\zeta\dd y.
\]
The second term has the same bound. The difference of the diffusion coefficients
contributes at most
\[
\norm{\delta v}_2\norm{\mathcal L\psi_m}_3\norm{\delta\psi}_6
 \le\eta_0\norm{\nabla\delta\psi}_2^2
       +C_{\eta_0}\norm{\nabla^2\psi_m}_3^2\norm{\delta\Phi}_2^2.
\]
The same estimate applies to $\kappa\delta v\Delta\zeta_m$.
Interpolation and the uniform second-order bound give
\[
 \norm{\nabla^2\psi_m}_3^2
 \le C\norm{\nabla^2\psi_m}_2\norm{\nabla^3\psi_m}_2
 \le C\bigl(1+D_2(Z_m)\bigr).
\]
For the transport difference, Hölder's inequality, Sobolev
embedding and Young's inequality give
\[
 \norm{\delta\psi}_6\norm{\nabla Z_m}_3\norm{\delta Z}_2
 \le\eta_0\norm{\nabla\delta\psi}_2^2
          +C_{\eta_0}\norm{\nabla Z_m}_3^2\norm{\delta Z}_2^2.
\]
For the remaining acoustic terms,
\[
\begin{aligned}
 &\left|\int\delta\zeta\nabla\Phi_m\cdot\delta\psi\dd y\right|
       +\left|\int(\delta\zeta)^2\Div\psi_m\dd y\right|\\
 &\quad\le C\norm{\delta\zeta}_6
                        \norm{\nabla Z_m}_3\norm{\delta Z}_2\\
 &\quad\le\eta_0\norm{\nabla\delta\zeta}_2^2
       +C_{\eta_0}\norm{\nabla Z_m}_3^2\norm{\delta Z}_2^2.
\end{aligned}
\]
For the acoustic term with variable temperature coefficient,
integration by parts gives
\[
\begin{aligned}
 &\int(\theta_n-\theta_*)\nabla\delta\Phi\cdot\delta\psi\dd y\\
 &\quad=-\int(\theta_n-\theta_*)\delta\Phi\Div\delta\psi\dd y
        -\int\delta\Phi\,\delta\psi\cdot
                          (\bar\theta'e_1+\nabla\zeta_n)\dd y.
\end{aligned}
\]
Using the $L^\infty$ bound for $\bar\theta'$ and the $L^3$
bound for $\nabla\zeta_n$, we obtain
\[
\begin{aligned}
 &\left|\int(\theta_n-\theta_*)
                  \nabla\delta\Phi\cdot\delta\psi\dd y\right|\\
 &\quad\le\eta_0\norm{\nabla\delta\psi}_2^2
       +C_{\eta_0,\delta}\bigl(1+\norm{\nabla\zeta_n}_3^2\bigr)
                                          \norm{\delta Z}_2^2.
\end{aligned}
\]
The remaining acoustic term containing the temperature difference satisfies
\[
\begin{aligned}
 &\left|\int\frac{\theta_n-\theta_*}{\theta_*}
                       \delta\zeta\Div\delta\psi\dd y\right|\\
 &\quad\le C\norm{\theta_n-\theta_*}_\infty
                         \norm{\delta\zeta}_2\norm{\nabla\delta\psi}_2\\
 &\quad\le\eta_0\norm{\nabla\delta\psi}_2^2
                      +C_{\eta_0,\delta}\norm{\delta\zeta}_2^2.
\end{aligned}
\]
By linearity of $\mathcal T$, the difference of the quadratic heating terms is
\[
\begin{aligned}
 &v_n\mathcal T(\psi_n):\nabla\psi_n
                   -v_m\mathcal T(\psi_m):\nabla\psi_m\\
 &\quad=v_n\mathcal T(\delta\psi):\nabla\psi_n
       +v_n\mathcal T(\psi_m):\nabla\delta\psi
       +\delta v\mathcal T(\psi_m):\nabla\psi_m.
\end{aligned}
\]
The first two terms, paired with $\delta\zeta/\theta_*$, are bounded by
\[
C\bigl(\norm{\nabla\psi_n}_3+\norm{\nabla\psi_m}_3\bigr)
        \norm{\nabla\delta\psi}_2\norm{\delta\zeta}_6,
\]
which is absorbed by the dissipation when the $H^2$ ball
is sufficiently small.
The term arising from the coefficient difference satisfies
\[
\begin{aligned}
 C\norm{\delta\Phi}_2\norm{\nabla\psi_m}_6^2
                       \norm{\delta\zeta}_6
 &\le\eta_0\norm{\nabla\delta\zeta}_2^2\\
 &\quad+C_{\eta_0}\norm{\nabla\psi_m}_6^4
                         \norm{\delta\Phi}_2^2.
\end{aligned}
\]
We also have
$\norm{\nabla\psi_m}_6\le C\norm{Z_m}_{H^2}$.
The profile contribution to the heating term has difference
\[
\begin{aligned}
 &v_n\mathcal T(\widetilde u):\nabla\psi_n
                   -v_m\mathcal T(\widetilde u):\nabla\psi_m\\
 &\quad=v_n\mathcal T(\widetilde u):\nabla\delta\psi
       +\delta v\mathcal T(\widetilde u):\nabla\psi_m.
\end{aligned}
\]
Since
$\norm{\mathcal T(\widetilde u)}_\infty\le C\delta^2$,
we have
\[
\begin{aligned}
 &\left|\int\frac{2\delta\zeta}{\theta_*}
       \bigl(v_n\mathcal T(\widetilde u):\nabla\delta\psi
          +\delta v\mathcal T(\widetilde u):\nabla\psi_m\bigr)\dd y\right|\\
 &\quad\le C\delta^2\left(
       \norm{\nabla\delta\psi}_2\norm{\delta\zeta}_2
       +\norm{\delta\Phi}_2\norm{\nabla\psi_m}_3
                                      \norm{\delta\zeta}_6\right)\\
 &\quad\le\eta_0\norm{\nabla(\delta\psi,\delta\zeta)}_2^2
       +C_{\eta_0,\delta}\bigl(1+\norm{\nabla\psi_m}_3^2\bigr)
                                             \norm{\delta Z}_2^2.
\end{aligned}
\]
All other profile terms are linear in $\delta Z$ and $\delta v$
and have bounded coefficients. Their contribution is therefore
bounded by $C_\delta\norm{\delta Z}_2^2$.
The forcing term satisfies
\[
 \left|\int\left(\theta_*\delta\Phi(\delta R)_\Phi
       +\delta\psi\cdot(\delta R)_\psi
       +\frac{c_v}{\theta_*}\delta\zeta(\delta R)_\zeta\right)\dd y\right|
 \le C\norm{\delta Z}_2^2+C\norm{\delta R}_2^2.
\]
Choose $\eta_0>0$ sufficiently small relative to the
parabolic coercivity constant, and then take the $H^2$
radius small enough to absorb the terms containing
$\nabla\Phi_n$ and the quadratic heating. This yields
\begin{equation}\label{g2:local-difference}
\frac{\dd}{\dd\tau}S_0(\delta Z)
 +c\norm{\nabla(\delta\psi,\delta\zeta)}_2^2
 \le C_\delta\{1+D_2(Z_n)+D_2(Z_m)\}\norm{\delta Z}_2^2
        +C\norm{\delta R}_2^2.
\end{equation}

By \eqref{g2:local-bound}, the coefficient on the right of
\eqref{g2:local-difference} has a time integral bounded
independently of the truncation indices. Gronwall's inequality gives
\[
\begin{aligned}
 &\sup_{0\le\tau\le T}\norm{\delta Z(\tau)}_2^2
   +\int_0^T\norm{\nabla(\delta\psi,\delta\zeta)}_2^2\dd\tau\\
 &\qquad\le C_{\delta,T}\left(
        \norm{(J_n-J_m)Z_0}_2^2
          +\norm{\delta R}_{L^2(0,T;L^2)}^2\right)
       \longrightarrow0.
\end{aligned}
\]
The interpolation estimate gives, for $0\le s<2$,
\[
 \norm{Z_n-Z_m}_{H^s}
 \le C\norm{Z_n-Z_m}_2^{1-s/2}
       \bigl(\norm{Z_n}_{H^2}+\norm{Z_m}_{H^2}\bigr)^{s/2}.
\]
Thus $Z_n$ converges strongly in $C([0,T];H^s)$.
Choose $s$ with $3/2<s<2$. Sobolev embedding gives convergence in
$C([0,T];L^\infty)$, and hence $v_n\to v$ in this space.
For example, the diffusion term satisfies
\[
 \norm{(v_n-v)\mathcal L\psi_n}_{L^2(0,T;L^2)}
 \le\norm{v_n-v}_{L^\infty((0,T)\times\R^3)}
                  \norm{\mathcal L\psi_n}_{L^2(0,T;L^2)}
 \longrightarrow0.
\]
The difference estimate also gives strong convergence of
$\nabla(\psi_n,\zeta_n)$ in $L^2(0,T;L^2)$. Consequently,
\[
\begin{aligned}
 &\norm{v_n\mathcal T(\psi_n):\nabla\psi_n
                    -v\mathcal T(\psi):\nabla\psi}_{L^1(0,T;L^1)}\\
 &\quad\le C\norm{v_n-v}_{L^\infty((0,T)\times\R^3)}
                       \norm{\nabla\psi_n}_{L^2(0,T;L^2)}^2\\
 &\qquad\quad+C_\delta
        \bigl(\norm{\nabla\psi_n}_{L^2(0,T;L^2)}
                  +\norm{\nabla\psi}_{L^2(0,T;L^2)}\bigr)
                  \norm{\nabla(\psi_n-\psi)}_{L^2(0,T;L^2)}
 \longrightarrow0.
\end{aligned}
\]
Passing to the limit in the equations gives a solution.
Weak compactness and \eqref{g2:local-bound} give the stated
bounds, while $R_n\to0$ removes the truncation error.
The difference estimate with zero forcing gives uniqueness.
The uniform $H^2$ bound and strong convergence in
$C([0,T];H^s)$ for every $s<2$ give weak continuity in $H^2$.
We also have
$(\psi_\tau,\zeta_\tau)\in L^2(0,T;H^1)$ and
$(\psi,\zeta)\in L^2(0,T;H^3)$. The Hilbert-triple theorem
then yields
$(\psi,\zeta)\in C([0,T];H^2)$. In particular,
\[
 \left|\frac{\dd}{\dd\tau}\norm{(\psi,\zeta)}_{H^2}^2\right|
 \le C\norm{(\psi_\tau,\zeta_\tau)}_{H^1}
                         \norm{(\psi,\zeta)}_{H^3}
 \in L^1(0,T).
\]
For the density component, the commutator norms of highest
order in the transport equation are bounded by
\[
 C\norm{\nabla\psi}_\infty\norm{\nabla^2\Phi}_2
 +C\norm{\nabla^2\psi}_6\norm{\nabla\Phi}_3
 +C_\delta\norm Z_{H^2}.
\]
The right-hand side is integrable in time. The source
$-\Div\psi-(\widetilde\rho'/\widetilde\rho)\psi_1$
also belongs to $L^1(0,T;H^2)$. Spatial mollification gives
the transport energy identity for each $|\alpha|\le2$:
\[
\begin{aligned}
 \frac12\frac{\dd}{\dd\tau}\norm{\partial^\alpha\Phi}_2^2
 &=\frac12\int\Div u\,|\partial^\alpha\Phi|^2\dd y\\
 &\quad-\int\partial^\alpha\Phi
       [\partial^\alpha,(u-\sigma e_1)\cdot\nabla]\Phi\dd y\\
 &\quad-\int\partial^\alpha\Phi\,
       \partial^\alpha\left(\Div\psi
                   +\frac{\widetilde\rho'}{\widetilde\rho}\psi_1\right)
                                                             \dd y.
\end{aligned}
\]
Summing over $|\alpha|\le2$ and using the smallness of
$\norm Z_{H^2}$, we obtain
\[
 \left|\frac{\dd}{\dd\tau}\norm\Phi_{H^2}^2\right|
 \le C_\delta\left(\norm Z_{H^2}^2
       +\norm Z_{H^2}\norm{\nabla\psi}_{H^2}\right)\in L^1(0,T).
\]
Absolute continuity of $\norm{\Phi}_{H^2}^2$, together with weak
continuity in $H^2$, gives $\Phi\in C([0,T];H^2)$.
The continuity equation then gives $\Phi_\tau\in C([0,T];H^1)$. If the $H^2$ norm is sufficiently
small, Sobolev embedding applied to
\[
 \rho=\widetilde\rho e^\Phi,\qquad
 \theta=\bar\theta+\zeta
\]
gives positive lower bounds for the density and temperature.

The local existence time depends only on the fixed profile and
the initial $H^2$ bound. Thus the solution can be restarted on
an interval of the same length whenever $\norm Z_{H^2}\le b$.

To propagate the negative-order regularity, we subtract the stationary
equations from the conservative system. The pressure, momentum
and energy flux differences then read
\[
\begin{aligned}
 p-\widetilde\rho\bar\theta
     &=\rho\zeta+\bar\theta r,\\
 \rho u\otimes u-\widetilde\rho\widetilde u\otimes\widetilde u
     &=\rho\psi\otimes u+\rho\widetilde u\otimes\psi
                              +r\widetilde u\otimes\widetilde u,\\
 (\rho E+p)u-(\widetilde\rho\widetilde E
                              +\widetilde\rho\bar\theta)\widetilde u
     &=(f_E+\rho\zeta+\bar\theta r)u
       +(\widetilde\rho\widetilde E
                              +\widetilde\rho\bar\theta)\psi.
\end{aligned}
\]
By the positive-state bounds and \eqref{g2:positive-equivalence},
each hyperbolic flux difference has $L^2$ norm at most
$C_\delta\norm f_2$. For the viscous energy flux, linearity of $\mathcal T$ gives
\[
\mathcal T(u)u-\mathcal T(\widetilde u)\widetilde u
 =\mathcal T(\psi)u+\mathcal T(\widetilde u)\psi.
\]
The bounds for $u$ and $\nabla\widetilde u$ therefore give
\[
 \norm{\mathcal T(\psi)u+\mathcal T(\widetilde u)\psi}_2
 \le C_\delta\bigl(\norm{\nabla\psi}_2+\norm\psi_2\bigr).
\]
The remaining viscous flux differences are $\mathcal T(\psi)$
and $\kappa\nabla\zeta$ in the momentum and energy equations,
respectively. Testing the conservative difference equations
against $\varphi\in C_c^\infty(\R^3;\R^5)$ yields
\[
 |\langle f_\tau,\varphi\rangle|
 \le C_\delta\left(\norm f_2
             +\norm{\nabla(\psi,\zeta)}_2\right)
                                      \norm{\nabla\varphi}_2.
\]
Taking the supremum in the definition
\eqref{eq:negative-norm}, we conclude
\[
 \norm{f_\tau}_{\dot H^{-1}}
 \le C_\delta\left(\norm f_2
                     +\norm{\nabla(\psi,\zeta)}_2\right).
\]
Both terms on the right are square-integrable in time, so
$f_\tau\in L^2(0,T;\dot H^{-1})$.
Since $f_0\in\dot H^{-1}$ and
$f_\tau\in L^2(0,T;\dot H^{-1})$, the fundamental theorem
of calculus in $\dot H^{-1}$ gives
\[
\begin{aligned}
 \norm{f(\tau)}_{\dot H^{-1}}
 &\le\norm{f_0}_{\dot H^{-1}}
              +\int_0^\tau\norm{f_s}_{\dot H^{-1}}\dd s\\
 &\le\norm{f_0}_{\dot H^{-1}}
              +T^{1/2}\norm{f_\tau}_{L^2(0,T;\dot H^{-1})}.
\end{aligned}
\]
This proves that $f\in H^1(0,T;\dot H^{-1})$.
The mass equation gives
\[
 r_\tau=-\Div\bigl(\rho u-\widetilde\rho\widetilde u
                            -\sigma r e_1\bigr)
                    \in L^2(0,T;L^2),
\]
and, since $\rho u-\widetilde\rho\widetilde u=f_m$,
\[
 \norm{r_\tau}_2
 \le C_\delta\norm f_{H^1}
 \le C_\delta\norm Z_{H^1}.
\]
Therefore,
$g_\tau=f_\tau+(0,d\nabla r_\tau,0)^{\mathsf T}$
belongs to $L^2(0,T;\dot H^{-1})$.
More precisely,
\[
 \norm g_{\dot H^{-1}}\le\norm f_{\dot H^{-1}}+C\norm r_2,
 \qquad
 \norm{g_\tau}_{\dot H^{-1}}
       \le\norm{f_\tau}_{\dot H^{-1}}+C\norm{r_\tau}_2.
\]
The two bounds together give
$g\in H^1(0,T;\dot H^{-1})$.
\end{proof}

\subsection{Proof of Theorem~\ref{g2:stability}}
\label{sec:stability-proof}

\begin{proposition}[A priori estimates]\label{g2:global-unit}
Let $S>0$ and let $Z$ be a local solution of
\eqref{ft:primitive} on $[0,S]$ in the class of
Lemma~\ref{g2:local}, with $f_0\in\dot H^{-1}$.
There exist constants $\delta_*,N_*>0$ and $C_\delta>0$,
independent of $S$, such that if
\[
 0<\delta\le\delta_*,\qquad
 \sup_{0\le\tau\le S}\norm{Z(\tau)}_{H^2}\le N_*,
\]
then
\begin{equation}\label{g2:unit-energy}
\begin{aligned}
&\sup_{0\le\tau\le S}
 \bigl(\norm{f(\tau)}_{\dot H^{-1}}^2+\norm{Z(\tau)}_{H^2}^2\bigr)\\
&\quad+\int_0^S\left\{\norm Z_{H^2}^2
                    +\norm{\nabla(\psi,\zeta)}_{H^2}^2\right\}\dd\tau
 \le C_\delta\bigl(\norm{f_0}_{\dot H^{-1}}^2
                                      +\norm{f_0}_{H^2}^2\bigr),
\end{aligned}
\end{equation}
\end{proposition}

The proof of Proposition~\ref{g2:global-unit} is given
in Section~\ref{sec:H2}.

\begin{proof}[Proof of Theorem~\ref{g2:stability}]
Work in the normalized stationary shock frame of
Section~\ref{sec:reformulation}. The temperature normalization
leaves the conservative variables unchanged. By Sobolev
embedding, the initial conservative perturbation is small in
$L^\infty$, so the initial density and temperature lie in
the fixed positive state neighborhood for sufficiently small
$a_{*,2}$. The norm equivalence
\eqref{g2:positive-equivalence}, together with the boundedness
of the Galilean map on $H^2$ and $\dot H^{-1}$, then gives
\[
 \norm{f_0}_{\dot H^{-1}}^2+\norm{f_0}_{H^2}^2
       +\norm{Z_0}_{H^2}^2\le C_\delta a_{*,2}^2,
\]
where $f_0$ and $Z_0$ denote the perturbations in the
normalized frame.

Let $\delta_*,N_*$ and $C_\delta$ be the constants in
Proposition~\ref{g2:global-unit}. Decrease $N_*$, if necessary,
so that it lies below the continuation radius in
Lemma~\ref{g2:local} and
\[
 CN_*\le\min\left\{1,\frac12\inf_{z\in\R}\bar\theta(z)\right\},
\]
where $C$ is the Sobolev embedding constant for $\Phi$ and
$\zeta$. After enlarging $C_\delta\ge1$ to cover the norm
equivalences, choose $a_{*,2}$ so that
$C_\delta^2a_{*,2}^2\le N_*^2/16$.

The initial data therefore satisfy
\[
 \norm{Z_0}_{H^2}\le\frac{N_*}{4},\qquad
 C_\delta\bigl(\norm{f_0}_{\dot H^{-1}}^2+\norm{f_0}_{H^2}^2\bigr)
                                       \le\frac{N_*^2}{4}.
\]
Lemma~\ref{g2:local} preserves the negative-order regularity
throughout the lifespan of the local solution. Set
\[
 S_*=\sup\left\{S>0:\ Z\text{ exists on }[0,S],\
                 \sup_{0\le\tau\le S}\norm{Z(\tau)}_{H^2}<N_*\right\}.
\]
Continuity, together with the strict initial bound
$\norm{Z_0}_{H^2}<N_*$, gives $S_*>0$. For every $S<S_*$,
Proposition~\ref{g2:global-unit} gives
\[
 \sup_{0\le\tau\le S}\norm{Z(\tau)}_{H^2}^2
       \le C_\delta\bigl(\norm{f_0}_{\dot H^{-1}}^2
                           +\norm{f_0}_{H^2}^2\bigr)
       \le\frac{N_*^2}{4}.
\]
Hence $\norm{Z(\tau)}_{H^2}\le N_*/2$ for
$0\le\tau<S_*$. If $S_*<\infty$, Lemma~\ref{g2:local} gives a common existence
time $T_0>0$ for initial data with $H^2$ norm at most $N_*/2$. Choose
$\tau_0<S_*$ with $S_*-\tau_0<T_0/2$ and restart the solution
at $\tau_0$. By uniqueness, it agrees with the original solution
on $[\tau_0,S_*)$.
Lemma~\ref{g2:local} also propagates the negative-order
regularity because $f(\tau_0)\in\dot H^{-1}$.
This solution is defined beyond $S_*$, with
\[
 \tau_0+T_0>S_*+T_0/2,\qquad
 \norm{Z(S_*)}_{H^2}
       =\lim_{\tau\uparrow S_*}\norm{Z(\tau)}_{H^2}\le N_*/2.
\]
Continuity in $H^2$ keeps the norm strictly below $N_*$ for a
short time beyond $S_*$. This contradicts the definition of
$S_*$, and hence $S_*=\infty$.
Letting $S\to\infty$ in \eqref{g2:unit-energy} and
using monotone convergence gives the global energy and dissipation estimate.
Lemma~\ref{g2:local} also gives
$f\in C([0,\infty);\dot H^{-1})$. The identities
$\rho=\widetilde\rho e^\Phi$ and
$\theta=\bar\theta+\zeta$, together with Sobolev embedding,
give
\[
 \widetilde\rho e^{-CN_*}\le\rho\le\widetilde\rho e^{CN_*},
 \qquad \bar\theta-CN_*\le\theta\le\bar\theta+CN_*.
\]
By the choice of $N_*$, we have
\[
 e^{-1}\inf_\R\widetilde\rho\le\rho\le e\sup_\R\widetilde\rho,
 \qquad
 \frac12\inf_\R\bar\theta\le\theta\le\frac32\sup_\R\bar\theta.
\]
These bounds hold uniformly in time.

Lemma~\ref{g2:local} gives local absolute continuity of $S_2$.
For the diffusion terms, integration by parts in the energy
identity of Lemma~\ref{ft:H2-lemma} yields
\[
 |S_2'(\tau)|\le C_\delta\bigl(D_2(\tau)+N(\tau)^2\bigr)
 \qquad\text{for almost every }\tau>0.
\]
Using the equivalence $S_2\sim N^2$, the global estimate
\eqref{g2:unit-energy} and the definition of $D_2$, we obtain
\[
 \int_0^\infty\bigl(S_2(\tau)+|S_2'(\tau)|\bigr)\dd\tau
 \le C_\delta\int_0^\infty
              \bigl(N(\tau)^2+D_2(\tau)\bigr)\dd\tau<\infty.
\]
Since $S_2'\in L^1(0,\infty)$, $S_2$ has a finite limit as
$\tau\to\infty$. This limit is zero because $S_2\ge0$ and
$S_2\in L^1(0,\infty)$.
The norm equivalence \eqref{g2:positive-equivalence} therefore gives
\mbox{$\norm{Z(\tau)}_{H^2}+\norm{f(\tau)}_{H^2}\to0$}.

We return to the original variables by undoing the temperature
normalization, the Galilean transformation and the change to the
stationary shock frame. The norm equivalences transfer
\eqref{g2:unit-energy} to \eqref{g2:stability-energy}.
Uniqueness and continuity in the original variables follow because
the transformations are invertible, spatial translations preserve
the norms, and time is unchanged.
\end{proof}

\subsection{Proof of Corollary~\ref{g2:theorem}}
\label{sec:limit}

\begin{proof}
\emph{Step 1. Scaling and global existence.}
In the physical variables of
\eqref{eq:NS}--\eqref{eq:coefficients}, introduce the scaling
\begin{equation}\label{eq:dilation}
 \tau=t/\eps,\qquad y=x/\eps,\qquad
 \widehat U(\tau,y)=U^\eps(\eps\tau,\eps y).
\end{equation}
The hats denote rescaled variables. By the chain rule,
\[
 \partial_\tau\widehat U=\eps\partial_tU^\eps,\qquad
 \nabla_y\widehat U=\eps\nabla_xU^\eps.
\]
The unscaled quantities on the left of the following identities
are evaluated at $(\eps\tau,\eps y)$:
\[
\begin{aligned}
 \mathcal T(u^\eps)&=\eps^{-1}\mathcal T(\widehat u),\\
 \eps\Div_x\mathcal T(u^\eps)
       &=\eps^{-1}\Div_y\mathcal T(\widehat u),\\
 \eps\Div_x\bigl(\mathcal T(u^\eps)u^\eps
                         +\kappa\nabla_x\theta^\eps\bigr)
       &=\eps^{-1}\Div_y\bigl(
          \mathcal T(\widehat u)\widehat u
                         +\kappa\nabla_y\widehat\theta\bigr).
\end{aligned}
\]
Inserting these identities into the equations, we cancel the
common factor $\eps^{-1}$ to obtain \eqref{eq:NS} with $\eps=1$.
The rescaled profile is
\[
 \widetilde U^\eps(\eps\tau,\eps y)
       =\widetilde U(y_1-\sigma\tau),
\]
with the same end states and shock speed. The corresponding
rescaled initial perturbation is
$\widehat w_0(y)=w_0^\eps(\eps y)$.
For $|\alpha|=k$,
\begin{equation}\label{eq:positive-scaling}
 \norm{\partial_y^\alpha\widehat w_0}_2^2
       =\eps^{2k-3}\norm{\partial_x^\alpha w_0^\eps}_2^2.
\end{equation}
For the negative-order norm, the scaling relation following
\eqref{g2:preparation} yields
\begin{equation}\label{eq:negative-scaling}
 \norm{\widehat w_0}_{\dot H^{-1}}^2
       =\eps^{-5}\norm{w_0^\eps}_{\dot H^{-1}}^2.
\end{equation}
Combining \eqref{eq:positive-scaling} and
\eqref{eq:negative-scaling}, we obtain
\begin{equation}\label{eq:initial-identity}
 \norm{\widehat w_0}_{H^2}^2
       +\norm{\widehat w_0}_{\dot H^{-1}}^2=A_{\eps,2}^2.
\end{equation}
Thus $A_{\eps,2}\le a_{*,2}$ is exactly the initial
smallness condition required by Theorem~\ref{g2:stability}
for the rescaled problem. The theorem provides a unique global
solution $\widehat U$. Rescaling gives a global solution $U^\eps$
with the same positive-state bounds and, for each fixed $\eps$,
the Sobolev regularity and continuity of $\widehat U$.

\medskip\noindent
\emph{Step 2. Uniform estimates in the physical variables.}

Applying \eqref{g2:stability-energy} to $\widehat U$
and using \eqref{eq:initial-identity}, we obtain
\[
\begin{aligned}
 &\sup_{\tau\ge0}\left(
       \norm{\widehat w(\tau)}_{\dot H^{-1}}^2
                        +\norm{\widehat Z(\tau)}_{H^2}^2\right)\\
 &\quad+\int_0^\infty\left(
       \norm{\widehat Z}_{H^2}^2
       +\norm{\nabla_y(\widehat\psi,\widehat\zeta)}_{H^2}^2\right)\dd\tau
 \le C_\delta A_{\eps,2}^2.
\end{aligned}
\]
At any fixed time,
\eqref{eq:positive-scaling}--\eqref{eq:negative-scaling}
apply to the perturbations. Since spatial translations preserve the
norms, the relation $t=\eps\tau$ gives, for $k=0,1,2$,
\[
 \eps^3\norm{\widehat w(\tau)}_{\dot H^{-1}}^2
       =\eps^{-2}\norm{w^\eps(t)}_{\dot H^{-1}}^2,
 \qquad
 \eps^3\norm{\nabla_y^k\widehat Z(\tau)}_2^2
       =\eps^{2k}\norm{\nabla_x^kZ^\eps(t)}_2^2.
\]
In the time integrals, the change of variables
$\dd\tau=\eps^{-1}\dd t$ gives
\begin{align*}
 \eps^3\int_0^\infty\norm{\nabla_y^k\widehat Z}_2^2\dd\tau
   &=\eps^{2k-1}\int_0^\infty\norm{\nabla_x^k Z^\eps}_2^2\dd t,\\
 \eps^3\int_0^\infty
       \norm{\nabla_y^{k+1}(\widehat\psi,\widehat\zeta)}_2^2\dd\tau
   &=\eps^{2k+1}\int_0^\infty
       \norm{\nabla_x^{k+1}(\psi^\eps,\zeta^\eps)}_2^2\dd t.
\end{align*}
Multiplying the rescaled estimate by $\eps^3$ and substituting these
identities, we obtain \eqref{g2:physical-energy}; here $t=\eps\tau$
identifies the two time suprema.

By \eqref{eq:pointwise-equivalence} and the Jacobian
$\dd x=\eps^3\dd y$, we have
\[
\begin{aligned}
 \norm{q^\eps(t)-\widetilde q^\eps(t)}_2^2
 &\le C\norm{Z^\eps(t)}_2^2
   =C\eps^3\norm{\widehat Z(\tau)}_2^2\\
 &\le C_\delta\eps^3 A_{\eps,2}^2.
\end{aligned}
\]
Applying the pointwise equivalence and Sobolev embedding in
the rescaled variables also yields
\[
 \norm{q^\eps(t)-\widetilde q^\eps(t)}_\infty
 \le C\norm{\widehat Z(\tau)}_\infty
 \le C\norm{\widehat Z(\tau)}_{H^2}
 \le C_\delta A_{\eps,2}.
\]
Taking the supremum over time, we obtain
\begin{align}
 \sup_{t\ge0}\norm{q^\eps(t)-\widetilde q^\eps(t)}_2
       &\le C_\delta\eps^{3/2}A_{\eps,2},
       \label{eq:L2-perturbation}\\
 \sup_{t\ge0}\norm{q^\eps(t)-\widetilde q^\eps(t)}_\infty
       &\le C_\delta A_{\eps,2}.
       \label{eq:Linfty-perturbation}
\end{align}
\Needspace{3\baselineskip}
\medskip\noindent
\emph{Step 3. Convergence to the Euler shock.}
For the profiles in \eqref{eq:profiles},
$0<M_\delta<\infty$ follows from \eqref{eq:layer-constant}.
The layer is independent of $x'$, so Fubini's theorem and the
change of variables $x_1=\sigma t+\eps y_1$ give
\[
 \int_\R\left|\widetilde q\left(\frac{x_1-\sigma t}{\eps}\right)
                    -q^s(t,x_1)\right|^2\dd x_1
 =\eps\int_\R|\widetilde q(y_1)-q^s(0,y_1)|^2\dd y_1
 =\eps M_\delta^2.
\]
Hence, for $0<|K|<\infty$,
\begin{equation}\label{eq:layer-identity}
 \norm{\widetilde q^\eps(t)-q^s(t)}_{L^2(\R\times K)}^2
           =\eps|K|M_\delta^2.
\end{equation}
Restricting to the cylinder and applying the reverse triangle
inequality, we obtain
\[
\begin{aligned}
&\left|\norm{q^\eps(t)-q^s(t)}_{L^2(\R\times K)}
       -\norm{\widetilde q^\eps(t)-q^s(t)}_{L^2(\R\times K)}\right|\\
&\qquad\le\norm{q^\eps(t)-\widetilde q^\eps(t)}_{L^2(\R\times K)}
       \le\norm{q^\eps(t)-\widetilde q^\eps(t)}_2.
\end{aligned}
\]
Together with \eqref{eq:L2-perturbation} and
\eqref{eq:layer-identity}, this proves \eqref{g2:cylinder}.
With $K$ and the shock fixed, $A_{\eps,2}\le a_{*,2}$ yields
\[
 \sup_{t\ge0}\left|
 \eps^{-1/2}\norm{q^\eps(t)-q^s(t)}_{L^2(\R\times K)}
                -|K|^{1/2}M_\delta\right|
 \le C_\delta\eps A_{\eps,2}\longrightarrow0.
\]
In particular, once
$C_\delta\eps a_{*,2}\le |K|^{1/2}M_\delta/2$,
\[
 \frac12|K|^{1/2}M_\delta\sqrt\eps
 \le\norm{q^\eps(t)-q^s(t)}_{L^2(\R\times K)}
 \le\frac32|K|^{1/2}M_\delta\sqrt\eps,
 \qquad t\ge0.
\]
This proves sharpness for each fixed shock and each such $K$.

On the region $|x_1-\sigma t|\ge\ell$,
\eqref{g2:profile-bounds} gives
\[
 |\widetilde q^\eps(t,x)-q^s(t,x)|
              \le C_\delta e^{-c_\delta\ell/\eps}.
\]
Combining the triangle inequality with \eqref{eq:Linfty-perturbation}
gives
\[
\begin{aligned}
 |q^\eps(t,x)-q^s(t,x)|
 &\le |q^\eps(t,x)-\widetilde q^\eps(t,x)|
                     +|\widetilde q^\eps(t,x)-q^s(t,x)|\\
 &\le C_\delta A_{\eps,2}
                     +C_\delta e^{-c_\delta\ell/\eps}.
\end{aligned}
\]
Taking the supremum over all $t\ge0$ and all points satisfying
$|x_1-\sigma t|\ge\ell$ gives \eqref{g2:away}. For each fixed $\ell>0$, the right-hand
side tends to zero when $A_{\eps,2}\to0$.
Proposition~\ref{prop:entropy}, using the estimates above,
gives the limit of the entropy-production measures and the
convergence of the conservative weak formulations.
\end{proof}

\subsection{Entropy production and the Euler limit}
\label{sec:entropy}

In the physical variables of
\eqref{eq:NS}--\eqref{eq:coefficients}, we use
\[
S(v,\theta)=\frac{R}{\gamma-1}\log\theta+R\log v,
\qquad [S]=S(v_+,\theta_+)-S(v_-,\theta_-),
\]
and write $S^\eps=S(v^\eps,\theta^\eps)$ and
$S^s=S(v^s,\theta^s)$. The mass Rankine--Hugoniot relation
shows that the mass flux in the shock frame is the same
at both end states. We denote it by
$j=\rho_\pm(u_{1\pm}-\sigma)$.

\begin{proposition}\label{prop:entropy}
Under the hypotheses of Corollary~\ref{g2:theorem},
the entropy-production measure is nonnegative and locally
finite on $(0,T)\times\R^3$ for every $T>0$. It is given by
\begin{equation}\label{eq:production}
\mathfrak M^\eps\dd x\dd t
 =\eps\left\{
 \frac{\mathcal T(u^\eps):\nabla u^\eps}{\theta^\eps}
 +\kappa\frac{|\nabla\theta^\eps|^2}{(\theta^\eps)^2}
          \right\}\dd x\dd t.
\end{equation}
As $\eps\to0$, these measures converge locally weakly to
\begin{equation}\label{eq:production-limit}
j[S]\,\delta_{x_1=\sigma t}\dd t\dd x'.
\end{equation}
More explicitly, the limit acts on a test function $\varphi$ by
$j[S]\int_0^T\int_{\R^2}\varphi(t,\sigma t,x')\dd x'\dd t$.
The solutions also converge to the stated shock in the
distributional mass, momentum, and total-energy equations.
\end{proposition}

\begin{proof}
For smooth solutions of \eqref{eq:NS} with positive density
and temperature, the continuity and internal-energy equations take the form
\[
\begin{aligned}
&(\partial_t+u^\eps\cdot\nabla)v^\eps
                                  =v^\eps\Div u^\eps,\\
&\frac{R\rho^\eps}{\gamma-1}
       (\partial_t+u^\eps\cdot\nabla)\theta^\eps
       +p^\eps\Div u^\eps
 =\eps\kappa\Delta\theta^\eps
       +\eps\mathcal T(u^\eps):\nabla u^\eps.
\end{aligned}
\]
Using the chain rule and $p^\eps/\theta^\eps=R\rho^\eps$,
we obtain
\[
\begin{aligned}
\rho^\eps(\partial_t+u^\eps\cdot\nabla)S^\eps
 &=\frac{R\rho^\eps}{(\gamma-1)\theta^\eps}
       (\partial_t+u^\eps\cdot\nabla)\theta^\eps
       +R\rho^\eps\Div u^\eps\\
 &=\eps\frac{\kappa\Delta\theta^\eps
                +\mathcal T(u^\eps):\nabla u^\eps}{\theta^\eps}.
\end{aligned}
\]
For the heat conduction term, we use the identity
\[
 \frac{\Delta\theta^\eps}{\theta^\eps}
 =\Div\left(\frac{\nabla\theta^\eps}{\theta^\eps}\right)
          +\frac{|\nabla\theta^\eps|^2}{(\theta^\eps)^2}.
\]
Expanding the conservative entropy derivatives and using the
continuity equation gives
\[
\begin{aligned}
 &\partial_t(\rho^\eps S^\eps)
                  +\Div(\rho^\eps u^\eps S^\eps)\\
 &\quad=\rho^\eps(\partial_t+u^\eps\cdot\nabla)S^\eps
       +S^\eps\bigl\{\partial_t\rho^\eps
                                  +\Div(\rho^\eps u^\eps)\bigr\}\\
 &\quad=\rho^\eps(\partial_t+u^\eps\cdot\nabla)S^\eps.
\end{aligned}
\]
Substituting the preceding identities gives
\begin{equation}\label{eq:entropy-balance}
\partial_t(\rho^\eps S^\eps)
 +\Div(\rho^\eps u^\eps S^\eps)
 =\eps\kappa\Div\left(\frac{\nabla\theta^\eps}{\theta^\eps}\right)
  +\mathfrak M^\eps.
\end{equation}
Fix $\eps>0$ and a bounded open set $B\subset\R^3$.
Apply Lemma~\ref{g2:local} in the rescaled variables and then
restore the smooth profile. This gives
\[
 (\rho^\eps,u^\eps,\theta^\eps)\in C([0,T];H^2(B)),\qquad
 (u^\eps,\theta^\eps)\in L^2(0,T;H^3(B)).
\]
The local bounds above justify \eqref{eq:entropy-balance}
by spatial mollification and the product estimates in
Lemma~\ref{ft:H2-lemma}.
The uniform positive-state bounds keep the entropy and its
derivatives bounded on the range of the solution.
The viscosity condition in \eqref{eq:coefficients} gives
\[
\mathcal T(u):\nabla u
 =2\mu\left|\frac{\nabla u+\nabla u^{\mathsf T}}2
                    -\frac{\Div u}{3}I\right|^2
  +\left(\lambda+\frac{2\mu}{3}\right)|\Div u|^2\ge0.
\]
Fix $0<T<\infty$ and a measurable set $K\subset\R^2$
with finite area. The change of variables
$x_1=\sigma t+\eps z$ gives
\[
\begin{aligned}
&\eps\int_0^T\int_{\R\times K}
 \bigl(|\nabla\widetilde u^\eps|^2
              +|\nabla\bar\theta^\eps|^2\bigr)\dd x\dd t\\
&\qquad=\eps\int_0^T\int_K\int_\R
        \eps^{-2}|(\widetilde u_1,\bar\theta)'(z)|^2
                                  \eps\,\dd z\dd x'\dd t\\
&\qquad=T|K|\norm{(\widetilde u_1,\bar\theta)'}_{L^2(\R)}^2.
\end{aligned}
\]
The profile bounds \eqref{g2:profile-bounds} imply
\[
 \norm{(\widetilde u_1,\bar\theta)'}_{L^2(\R)}^2
 \le C\delta^4\int_\R e^{-2c\delta|z|}\dd z
 \le C\delta^3.
\]
Write each gradient as the sum of its profile and perturbation
parts, as in \eqref{eq:perturbations}. The $k=0$
dissipation of velocity and temperature in
\eqref{g2:physical-energy} then gives
\begin{equation}\label{eq:local-dissipation}
\begin{aligned}
&\eps\int_0^T\int_{\R\times K}
       \bigl(|\nabla u^\eps|^2+|\nabla\theta^\eps|^2\bigr)\dd x\dd t\\
&\quad\le C\eps\int_0^T
             \norm{\nabla(\psi^\eps,\zeta^\eps)}_2^2\dd t
       +CT|K|\norm{(\widetilde u_1,\bar\theta)'}_{L^2(\R)}^2\\
&\quad\le C_\delta\bigl(\eps^3A_{\eps,2}^2+T|K|\bigr).
\end{aligned}
\end{equation}
The uniform state bounds imply
\[
 0\le\mathfrak M^\eps
   \le C_\delta\eps
           \bigl(|\nabla u^\eps|^2+|\nabla\theta^\eps|^2\bigr).
\]
By \eqref{eq:local-dissipation}, the measures in
\eqref{eq:production} have uniformly bounded mass on compact sets.
The heat flux satisfies
\[
\begin{aligned}
 &\int_0^T\int_{\R\times K}
       \left|\frac{\eps\kappa\nabla\theta^\eps}{\theta^\eps}\right|^2
                                                       \dd x\dd t\\
 &\qquad\le C_\delta\eps
       \left(\eps\int_0^T\int_{\R\times K}
                              |\nabla\theta^\eps|^2\dd x\dd t\right).
\end{aligned}
\]
Substituting the local dissipation bound gives
\begin{equation}\label{eq:heat-flux}
\norm{\frac{\eps\kappa\nabla\theta^\eps}{\theta^\eps}}
 _{L^2((0,T)\times(\R\times K))}^2
 \le C_\delta\bigl(\eps^4A_{\eps,2}^2+\eps T|K|\bigr)
 \longrightarrow0.
\end{equation}

On the state range under consideration, the entropy density
and entropy flux are Lipschitz functions of $q$.
Let $B\subset\R^3$ be bounded, and choose a bounded
$K\subset\R^2$ such that $B\subset\R\times K$.
By Cauchy--Schwarz, the triangle inequality and
\eqref{eq:layer-identity},
\[
\begin{aligned}
 \norm{q^\eps(t)-q^s(t)}_{L^1(B)}
 &\le |B|^{1/2}\norm{q^\eps(t)-q^s(t)}_{L^2(\R\times K)}\\
 &\le |B|^{1/2}\left(
       \norm{q^\eps(t)-\widetilde q^\eps(t)}_2
                         +|K|^{1/2}M_\delta\sqrt\eps\right).
\end{aligned}
\]
The Lipschitz bounds and \eqref{eq:L2-perturbation}
therefore give
\[
\begin{aligned}
\sup_{0\le t\le T}\bigl(&
 \norm{\rho^\eps S^\eps-\rho^sS^s}_{L^1(B)}
 +\norm{\rho^\eps u^\eps S^\eps-\rho^su^sS^s}_{L^1(B)}\bigr)\\
 &\le C_{\delta,B}\bigl(\sqrt\eps+\eps^{3/2}A_{\eps,2}\bigr)
 \longrightarrow0.
\end{aligned}
\]
Take $\varphi\in C_c^\infty((0,T)\times\R^3)$.
Testing \eqref{eq:entropy-balance} against $\varphi$ and
integrating by parts, we obtain
\[
\begin{aligned}
\int_0^T\int_{\R^3}\mathfrak M^\eps\varphi\dd x\dd t
 &=-\int_0^T\int_{\R^3}\rho^\eps S^\eps
                                 \partial_t\varphi\dd x\dd t\\
 &\quad-\int_0^T\int_{\R^3}\rho^\eps u^\eps S^\eps
                                 \cdot\nabla\varphi\dd x\dd t\\
 &\quad+\eps\kappa\int_0^T\int_{\R^3}
         \frac{\nabla\theta^\eps}{\theta^\eps}
                                 \cdot\nabla\varphi\dd x\dd t.
\end{aligned}
\]
The local $L^1$ convergence gives the limits of the first
two terms, while the last tends to zero by \eqref{eq:heat-flux}.
For the limiting shock, integration by parts over the two
half-spaces separated by $x_1=\sigma t$ gives
\[
\begin{aligned}
 &-\int_0^T\int_{\R^3}\bigl(
       \rho^sS^s\partial_t\varphi
             +\rho^su^sS^s\cdot\nabla\varphi\bigr)\dd x\dd t\\
 &\quad=\int_0^T\int_{\R^2}
       \bigl([\rho u_1S]-\sigma[\rho S]\bigr)
                      \varphi(t,\sigma t,x')\dd x'\dd t.
\end{aligned}
\]
Hence the limiting entropy balance is
\[
\partial_t(\rho^sS^s)+\partial_1(\rho^s u_1^sS^s)
 =\bigl([\rho u_1S]-\sigma[\rho S]\bigr)
                       \delta_{x_1=\sigma t}
 =j[S]\delta_{x_1=\sigma t}.
\]
For the last equality, the common mass flux gives
\[
\begin{aligned}
[\rho u_1S]-\sigma[\rho S]
 &=\rho_+(u_{1+}-\sigma)S(v_+,\theta_+)
    -\rho_-(u_{1-}-\sigma)S(v_-,\theta_-)\\
 &=j[S].
\end{aligned}
\]
For $\varphi\in C_c((0,T)\times\R^3)$,
choose $\varphi_n\in C_c^\infty((0,T)\times\R^3)$
that converge uniformly to $\varphi$. We may take all their
supports in one fixed cylinder $(0,T)\times B$, with $B$ bounded.
Since convergence holds for each $\varphi_n$, the uniform
local mass bound gives
\[
\begin{aligned}
 &\limsup_{\eps\to0}\left|
       \int_0^T\int_{\R^3}\mathfrak M^\eps\varphi\dd x\dd t
       -j[S]\int_0^T\int_{\R^2}
                  \varphi(t,\sigma t,x')\dd x'\dd t\right|\\
 &\qquad\le C_{\delta,T,B}\norm{\varphi-\varphi_n}_\infty
                         \longrightarrow0\qquad(n\to\infty).
\end{aligned}
\]
Letting $n\to\infty$ proves local weak convergence of the
full family of entropy-production measures.

For the profile with $\eps=1$, the constant mass flux gives
\[
\begin{aligned}
 &-\sigma\bigl(\widetilde\rho S(\widetilde v,\bar\theta)\bigr)'
       +\bigl(\widetilde\rho\widetilde u_1
                             S(\widetilde v,\bar\theta)\bigr)'\\
 &\quad=\bigl\{\widetilde\rho(\widetilde u_1-\sigma)
                             S(\widetilde v,\bar\theta)\bigr\}'
       =j\frac{\dd}{\dd z}S(\widetilde v,\bar\theta).
\end{aligned}
\]
The profile version of \eqref{eq:entropy-balance} is
\[
 \frac{\dd}{\dd z}\left\{
       jS(\widetilde v,\bar\theta)
                      -\kappa\frac{\bar\theta'}{\bar\theta}\right\}
 =\frac{\nu(\widetilde u_1')^2}{\bar\theta}
                      +\kappa\frac{(\bar\theta')^2}{\bar\theta^2}.
\]
The decay of the profile implies that
$\bar\theta'/\bar\theta$ vanishes at both end states.
Integration over $\R$ gives
\[
 j[S]=\int_\R\left\{
       \frac{\nu(\widetilde u_1')^2}{\bar\theta}
                      +\kappa\frac{(\bar\theta')^2}{\bar\theta^2}
                      \right\}\dd z\ge0.
\]

It remains to pass to the limit in the conservative equations.
For the Euler shock, the conservative variables are
\[
 \rho^s=(v^s)^{-1},\qquad p^s=R\rho^s\theta^s,\qquad
 E^s=\frac{R\theta^s}{\gamma-1}+\frac{|u^s|^2}{2}.
\]
The conservative variables and Euler fluxes are Lipschitz
functions of $q$ on the state range. For the energy flux,
we use the decomposition
\[
\begin{aligned}
 &(\rho^\eps E^\eps+p^\eps)u^\eps
                         -(\rho^sE^s+p^s)u^s\\
 &\quad=(\rho^\eps E^\eps-\rho^sE^s+p^\eps-p^s)u^\eps
                         +(\rho^sE^s+p^s)(u^\eps-u^s).
\end{aligned}
\]
The uniform state bounds show that this difference is bounded
pointwise by $C_\delta|q^\eps-q^s|$.
The mass and momentum fluxes admit analogous decompositions,
so the same bound gives their local $L^1$ convergence.
At $t=0$, \eqref{g2:cylinder} also gives local $L^1$
convergence of the conservative initial data to those of $q^s$.
Since $|\mathcal T(u^\eps)|\le C|\nabla u^\eps|$,
we obtain from \eqref{eq:local-dissipation}
\[
 \eps^2\int_0^T\int_{\R\times K}
 \bigl(|\mathcal T(u^\eps)|^2
                  +\kappa^2|\nabla\theta^\eps|^2\bigr)\dd x\dd t
 \le C_\delta\bigl(\eps^4A_{\eps,2}^2+\eps T|K|\bigr)
 \longrightarrow0.
\]
The uniform velocity bound also gives
$\eps\mathcal T(u^\eps)u^\eps\to0$ in local
space-time $L^2$.
Take test functions
$\varphi\in C_c^\infty([0,T)\times\R^3)$ and
$\xi\in C_c^\infty([0,T)\times\R^3;\R^3)$.
The weak formulations of \eqref{eq:NS} read
\begin{align*}
 0&=\int_{\R^3}\rho_0^\eps\varphi(0)\dd x\\*
  &\quad+\int_0^T\int_{\R^3}
       \bigl(\rho^\eps\partial_t\varphi
                    +\rho^\eps u^\eps\cdot\nabla\varphi\bigr)
                                                       \dd x\dd t,\\[2pt]
 0&=\int_{\R^3}\rho_0^\eps u_0^\eps\cdot\xi(0)\dd x\\*
  &\quad+\int_0^T\int_{\R^3}
       \bigl(\rho^\eps u^\eps\cdot\partial_t\xi
           +\rho^\eps u^\eps\otimes u^\eps:\nabla\xi
                                      +p^\eps\Div\xi\bigr)\dd x\dd t\\*
  &\quad-\eps\int_0^T\int_{\R^3}
                           \mathcal T(u^\eps):\nabla\xi\dd x\dd t,\\[2pt]
 0&=\int_{\R^3}\rho_0^\eps E_0^\eps\varphi(0)\dd x\\*
  &\quad+\int_0^T\int_{\R^3}
       \bigl(\rho^\eps E^\eps\partial_t\varphi
              +(\rho^\eps E^\eps+p^\eps)u^\eps\cdot\nabla\varphi\bigr)
                                                       \dd x\dd t\\*
  &\quad-\eps\int_0^T\int_{\R^3}
       \bigl(\mathcal T(u^\eps)u^\eps
                      +\kappa\nabla\theta^\eps\bigr)
                                  \cdot\nabla\varphi\dd x\dd t.
\end{align*}

The local $L^1$ convergence gives the limits of the initial-data
and Euler-flux terms. Choose a finite-area set $K\subset\R^2$
with both test-function supports in $\R\times K$.
By Cauchy--Schwarz,
\[
\begin{aligned}
 &\left|\eps\int_0^T\int_{\R^3}
                         \mathcal T(u^\eps):\nabla\xi\dd x\dd t\right|\\
 &\quad+\left|\eps\int_0^T\int_{\R^3}
       \bigl(\mathcal T(u^\eps)u^\eps+\kappa\nabla\theta^\eps\bigr)
                                  \cdot\nabla\varphi\dd x\dd t\right|\\
 &\qquad\le C_\delta\bigl(\eps^2A_{\eps,2}
                                 +\sqrt{\eps T|K|}\bigr)
       \left(\norm{\nabla\xi}_{L^2((0,T)\times\R^3)}
           +\norm{\nabla\varphi}_{L^2((0,T)\times\R^3)}\right)
       \longrightarrow0.
\end{aligned}
\]
The dissipative terms vanish as $\eps\to0$, so passing to
the limit in the weak formulations gives the Euler equations
with initial data $q^s$. The Rankine--Hugoniot relations
\eqref{eq:RH} then give
\[
\begin{aligned}
 \partial_t\rho^s+\partial_1(\rho^su_1^s)
   &=\bigl(-\sigma[\rho]+[\rho u_1]\bigr)
                                      \delta_{x_1=\sigma t}=0,\\
 \partial_t(\rho^su_1^s)
       +\partial_1\bigl(\rho^s(u_1^s)^2+p^s\bigr)
   &=\bigl(-\sigma[\rho u_1]+[\rho u_1^2+p]\bigr)
                                      \delta_{x_1=\sigma t}=0,\\
 \partial_t(\rho^sE^s)
       +\partial_1\bigl((\rho^sE^s+p^s)u_1^s\bigr)
   &=\bigl(-\sigma[\rho E]+[(\rho E+p)u_1]\bigr)
                                      \delta_{x_1=\sigma t}=0.
\end{aligned}
\]
The transverse momentum equations are also satisfied, since
$u_2^s=u_3^s=0$ and $p^s$ is independent of the transverse
variables.
\end{proof}

\section{Uniform-in-time \texorpdfstring{$H^2$}{H2} estimates}
\label{sec:H2}

Throughout this section, we work with the normalized unit-viscosity system
\eqref{ft:primitive} on an interval $[0,S]$ for which
$\sup_{0\le\tau\le S}N(\tau)\le N_*$.
All constants below are independent of $S$.

\subsection{Effective perturbation and elliptic energy}
\label{sec:negative}

We use the normalized variables of Section~\ref{sec:reformulation}.
Define the effective conservative perturbation
\begin{equation}\label{ft:effective}
d=\min\left\{\frac{\nu}{2\rho_*},\frac{\kappa}{\rho_*}\right\}>0,
\qquad g=f+(0,d\nabla r,0)^{\mathsf T}.
\end{equation}
For $\varphi\in C_c^\infty(\R^3;\R^3)$, integration
by parts gives
\[
 |\langle\nabla r,\varphi\rangle|
       =\left|\int r\Div\varphi\dd y\right|
       \le\sqrt3\norm r_2\norm{\nabla\varphi}_2.
\]
Taking the supremum over the test functions in
\eqref{eq:negative-norm} gives
\[
 \norm{\nabla r}_{\dot H^{-1}}\le C\norm r_2,
 \qquad
 \norm g_{\dot H^{-1}}\le\norm f_{\dot H^{-1}}+C\norm r_2.
\]
Thus $f\in H^2\cap\dot H^{-1}$ implies
$g\in L^2\cap\dot H^{-1}$. We also use
\[
 \norm g_2\le\norm f_2+d\norm{\nabla r}_2,
 \qquad
 \norm f_{\dot H^{-1}}\le\norm g_{\dot H^{-1}}+C\norm r_2.
\]
Define the entropy function
$\mathcal S(U)=-\rho(c_v\log\theta-\log\rho)$ and let
$M=(\nabla_U^2\mathcal S(\widetilde U))^{-1}$.
The inverse entropy Hessian in
\cite[equation (2.46)]{HuangXu2026} satisfies
\[
h^{\mathsf T}Mh=\widetilde\rho\left\{
 \bigl(h_\rho+\widetilde u_1h_{m,1}
             +(c_v\bar\theta+\widetilde u_1^2/2)h_E\bigr)^2
 +\bar\theta|h_m+\widetilde u_1h_Ee_1|^2
 +c_v\bar\theta^2h_E^2\right\}.
\]
The uniform upper and lower profile bounds give
\mbox{$c|h|^2\le h^{\mathsf T}Mh\le C|h|^2$}.
We fix the exponent $k>0$ as in
\cite[Proposition~3.1, Step~2]{HuangXu2026} and define
\[
\omega=(\widetilde\rho/\rho_*)^k,\qquad
Jh=h_m+\widetilde u_1h_Ee_1.
\]
The exponent $k$ is chosen independently of $\delta$,
after which we take the shock strength sufficiently small. With $k$ fixed,
$\omega$ is bounded above and below by positive constants
along the profile.
For $h,\varphi\in\dot H^1(\R^3;\R^5)$, define the bilinear form
\begin{equation}\label{ft:elliptic-form}
\begin{aligned}
\mathfrak a(h,\varphi)
 &=\int\omega\sum_{j=1}^3(\partial_jh)^{\mathsf T}M\partial_j\varphi\dd y\\
 &\quad+\int\omega\widetilde\rho\bar\theta
          \nabla\times(Jh)\cdot\nabla\times(J\varphi)\dd y,
\qquad \langle\mathcal A h,\varphi\rangle=\mathfrak a(h,\varphi).
\end{aligned}
\end{equation}
Since the profile depends only on $y_1$, we have
\[
\begin{aligned}
\nabla\times(\widetilde u_1h_Ee_1)
 &=(\widetilde u_1'h_Ee_1+\widetilde u_1\nabla h_E)\times e_1\\
 &=\widetilde u_1\nabla h_E\times e_1
   =\widetilde u_1(0,\partial_3h_E,-\partial_2h_E)^{\mathsf T}.
\end{aligned}
\]
Thus $|\nabla\times(Jh)|\le C|\nabla h|$, and the bounds
on the profile and the weight give
\[
 |\mathfrak a(h,\varphi)|\le C\norm{\nabla h}_2\norm{\nabla \varphi}_2,
 \qquad
 c\norm{\nabla h}_2^2\le\mathfrak a(h,h)\le C\norm{\nabla h}_2^2.
\]
With the normalization in Section~\ref{sec:reformulation},
$U_E=\rho E$ is the conservative energy variable used by
Huang--Xu. Our matrix and bilinear form therefore agree with
\cite[equations (2.46)--(2.47), Lemma~2.2]{HuangXu2026}.
These bounds and the Lax--Milgram theorem show that
$\mathcal A:\dot H^1\to\dot H^{-1}$ is an isomorphism.
Since the coefficients depend only on $y_1$, $\mathcal A$
is independent of $\tau$.
Integration by parts in \eqref{ft:elliptic-form} gives
\[
 \mathcal A h=-\sum_{j=1}^3\partial_j(\omega M\partial_jh)
 +J^{\mathsf T}\nabla\times
       \bigl(\omega\widetilde\rho\bar\theta
                         \nabla\times(Jh)\bigr)
\]
as an identity in distributions. Here $J^{\mathsf T}$ is
the transpose of $J$.
For $g\in\dot H^{-1}$, define
\begin{equation}\label{g2:negative-energy}
h=\mathcal A^{-1}g,\qquad
E_-=\frac12\langle g,h\rangle=\frac12\mathfrak a(h,h)
       \sim\norm g_{\dot H^{-1}}^2.
\end{equation}
The norm equivalence in \eqref{g2:negative-energy} follows from
\cite[Lemma~2.2, equation (2.48)]{HuangXu2026} under this normalization.
The Sobolev inequality also gives
$\norm h_6\le C\norm{\nabla h}_2\le C\norm g_{\dot H^{-1}}$.
We henceforth regard $h$ as its $L^6$ representative.

By Lemma~\ref{lm:profile}, the coefficients satisfy
\[
 \omega'=k\omega\frac{\widetilde\rho'}{\widetilde\rho},
 \qquad |\omega'|+|M'|+|J'|\le C_k|\widetilde u_1'|.
\]
With $k$ fixed and the shock sufficiently weak, the estimate in
\cite[Lemma~2.3, equation (2.49)]{HuangXu2026} applies to
$g\in L^2\cap\dot H^{-1}$. It gives $D^2h\in L^2$ and
\begin{equation}\label{g2:elliptic-bounds}
\begin{gathered}
D_0=\norm{D^2h}_2^2,\qquad
D_s=\int\omega|\widetilde u_1'|\,|\nabla h|^2\dd y,\\
\norm{D^2h}_2\le C\bigl(\norm g_2+
                       \norm{\widetilde u_1'\nabla h}_2\bigr),\qquad
\norm g_2^2\le C(D_0+\delta^2D_s),\qquad
\norm{\widetilde u_1'\nabla h}_2\le C\delta D_s^{1/2}.
\end{gathered}
\end{equation}
Here $D^2h$ denotes the full Hessian, with the sum taken over
ordered pairs of spatial indices. Both the Hessian estimate and the bound for $\norm g_2$
follow from the cited lemma under $g\in L^2\cap\dot H^{-1}$
and the profile coefficient bounds. No derivative of $g$ is used.
Expanding $\mathcal A$ gives
\[
 |g|\le C|D^2h|+C|\widetilde u_1'|\,|\nabla h|,
 \qquad
 \norm{\widetilde u_1'\nabla h}_2^2
       \le C\delta^2D_s.
\]
Applying the Sobolev inequality to $\nabla h$ gives
$\norm{\nabla h}_6\le CD_0^{1/2}$.

The primitive perturbations satisfy
\begin{equation}\label{g2:recovery-identities}
\begin{aligned}
r&=\widetilde\rho(e^\Phi-1),\qquad
 \nabla r=\rho\nabla\Phi+\widetilde\rho'(e^\Phi-1)e_1,\\
\rho\psi&=g_m-\widetilde u r-d\nabla r,\\
c_v\rho\zeta&=f_E-r\left(c_v\bar\theta+\tfrac12|\widetilde u|^2\right)
              -\rho\widetilde u\cdot\psi-\tfrac12\rho|\psi|^2.
\end{aligned}
\end{equation}
The argument in \cite[Lemma~3.2, equations (3.18)--(3.19)]{HuangXu2026}
requires only the positive-state and profile bounds, smallness of
$\norm Z_\infty$, the algebraic relations above with $r=g_\rho$,
and \eqref{g2:elliptic-bounds}. These hypotheses follow from the preceding bounds and
$H^2(\R^3)\hookrightarrow L^\infty$. The case $m=2$ of that argument
then yields
\begin{equation}\label{g2:recovery}
\begin{aligned}
\norm Z_{H^2}^2&\le C(D_0+D_2+\delta^2D_s),\\
\norm{\nabla r}_2^2&\le CD_2+C\delta^4(D_0+\delta^2D_s),\\
\norm{\nabla f}_2^2&\le CD_2+C\delta^4\norm Z_2^2.
\end{aligned}
\end{equation}

\subsection{Negative-order estimate}

\begin{lemma}\label{g2:negative-lemma}
Let $Z$ be a local solution from Lemma~\ref{g2:local}
with $f_0\in\dot H^{-1}$. Suppose that $\delta+N$ is
sufficiently small. Then
\begin{equation}\label{g2:negative-inequality}
E_-'+c(D_0+D_s)
 \le C\norm{\nabla r}_2^2
        +CN(D_0+D_2+\norm Z_2^2).
\end{equation}
The constants $c,C$ depend only on the fixed coefficients
and the positive state neighborhood.
\end{lemma}

\begin{proof}
Lemma~\ref{g2:local} gives $g\in H^1(0,T;\dot H^{-1})$.
Since $\mathcal A$ is time-independent and self-adjoint, the
preceding time regularity allows us to apply
\cite[equation (2.50)]{HuangXu2026} and obtain
\begin{equation}\label{g2:negative-identity}
E_-'=\langle g_\tau,h\rangle
\end{equation}
for almost every $\tau$. Density in $\dot H^1\cap\dot H^2$
justifies the integrations by parts below.

In the stationary shock frame, let $F_j$ denote the Euler
flux in the $j$-th direction and set
\[
 F_j(U)=\bigl(\rho u_j,\rho u u_j+pe_j,
                        (\rho E+p)u_j\bigr)^{\mathsf T},
 \qquad A_j=DF_j(\widetilde U)-\sigma\delta_{j1}I.
\]
Substituting $f_m=g_m-d\nabla r$ into the mass equation gives
\[
 r_\tau-\sigma\partial_1r+\Div g_m=d\Delta r.
\]
Differentiating the momentum component of \eqref{ft:effective}
and using this identity, we find
\[
\begin{aligned}
 \partial_\tau g_m-\sigma\partial_1g_m
 &=\partial_\tau f_m-\sigma\partial_1f_m
                         +d\nabla(r_\tau-\sigma\partial_1r)\\
 &=\partial_\tau f_m-\sigma\partial_1f_m
                         -d\nabla\Div g_m+d^2\nabla\Delta r.
\end{aligned}
\]
Differentiating the Euler fluxes for the correction in
\eqref{ft:effective} gives
\[
\begin{aligned}
&\sum_{j=1}^3\partial_j
       \{DF_j(\widetilde U)(0,d\nabla r,0)^{\mathsf T}\}\\
&\quad=d\begin{pmatrix}
 \Delta r\\
 \partial_1(\widetilde u_1\nabla r)
       +e_1\Div(\widetilde u_1\nabla r)
       -c_v^{-1}\nabla(\widetilde u_1\partial_1r)\\
 \Div\{((c_v+1)\bar\theta+\widetilde u_1^2/2)\nabla r\}
       -c_v^{-1}\partial_1(\widetilde u_1^2\partial_1r)
\end{pmatrix}.
\end{aligned}
\]
Combining the mass, momentum and energy equations gives
\begin{equation}\label{g2:effective-system}
g_\tau+\sum_{j=1}^3\partial_j(A_jg)=\mathcal B g+\mathcal N_g.
\end{equation}
The corresponding linearized velocity and temperature perturbations are
\[
\begin{aligned}
 \widehat u
 &=\frac{f_m-\widetilde u r}{\widetilde\rho}
  =\frac{g_m-\widetilde u_1re_1-d\nabla r}{\widetilde\rho},\\
 \widehat\theta
 &=\frac{f_E-\widetilde u\cdot f_m
                 +(\widetilde u_1^2/2-c_v\bar\theta)r}
                                      {c_v\widetilde\rho}\\
 &=\frac{g_E-\widetilde u_1g_{m,1}
 +(\widetilde u_1^2/2-c_v\bar\theta)r+d\widetilde u_1\partial_1r}
                                      {c_v\widetilde\rho}.
\end{aligned}
\]
The linear diffusion operator $\mathcal B$ is given by
\begin{equation}\label{g2:linear-diffusion}
\begin{aligned}
(\mathcal B g)_\rho={}&d\Delta r,\\
(\mathcal B g)_m={}&\mathcal L\widehat u-d\nabla\Div g_m+d^2\nabla\Delta r
 +d\partial_1(\widetilde u_1\nabla r)\\
 &+de_1\Div(\widetilde u_1\nabla r)
             -\frac d{c_v}\nabla(\widetilde u_1\partial_1r),\\
(\mathcal B g)_E={}&\kappa\Delta\widehat\theta
 +\Div\{\mathcal T(\widetilde u)\widehat u
                         +\mathcal T(\widehat u)\widetilde u\}\\
 &+d\Div\{((c_v+1)\bar\theta+\widetilde u_1^2/2)\nabla r\}
                 -\frac d{c_v}\partial_1(\widetilde u_1^2\partial_1r).
\end{aligned}
\end{equation}
Under the change of variables in Section~\ref{sec:reformulation},
the operator agrees with that in
\cite[equations (3.2)--(3.3)]{HuangXu2026}; the moving-frame
term is included in $A_1$. Pairing \eqref{g2:effective-system}
with $h$ and integrating the Euler fluxes by parts, we obtain
\[
 E_-'=\left\langle g,\sum_{j=1}^3A_j^{\mathsf T}\partial_jh\right\rangle
                  +\langle\mathcal B g,h\rangle
                  +\langle\mathcal N_g,h\rangle.
\]

For the Euler contribution, we use the entropy symmetry
$MA_j^{\mathsf T}=A_jM$ and the modified curl in $\mathcal A$.
After the normalization in Section~\ref{sec:reformulation},
the pairing is the same as that in
\cite[Proposition~3.1, Step~2, equations (3.7)--(3.11)]{HuangXu2026},
with the moving-frame term absorbed into $A_1$. The second derivatives of $h$ cancel, and the resulting pairing
is well defined for \mbox{$h\in\dot H^1\cap\dot H^2$} using only the
entropy symmetry and the profile bounds. Choosing $k$ as in the
cited proposition gives
\[
 \left\langle\mathcal A h,
       \sum_{j=1}^3A_j^{\mathsf T}\partial_jh\right\rangle
       \le-2c_E D_s+C_k\delta D_s.
\]
Choosing $\delta$ sufficiently small that $C_k\delta\le c_E$,
we obtain
\begin{equation}\label{g2:Euler-pairing}
\left\langle\mathcal A h,
                \sum_{j=1}^3A_j^{\mathsf T}\partial_jh\right\rangle
                                                    \le-c_E D_s.
\end{equation}

For the diffusion contribution,
\cite[Proposition~3.1, Steps~3--4]{HuangXu2026} applied to
\eqref{g2:linear-diffusion} gives
\begin{equation}\label{g2:diffusive-pairing}
\langle\mathcal B g,h\rangle
       \le-c_0D_0+\frac{c_E}{2}D_s+C\norm{\nabla r}_2^2.
\end{equation}
At this regularity, we integrate by parts to place the derivatives
on the elliptic potential. Its density--energy part at the frozen state
$(\rho_*,0,\theta_*)$ is
\begin{equation}\label{g2:scalar-square}
-\rho_*d\norm{D^2\{h_\rho+(c_v+\tfrac12)\theta_*h_E\}}_2^2
 -\theta_*^2\left(\kappa-\frac14\rho_*d\right)\norm{D^2h_E}_2^2.
\end{equation}
By the definition of $d$,
$\kappa-\rho_*d/4\ge3\kappa/4>0$. The momentum part takes the form
\[
 -\rho_*\theta_*\left\{
 \left(\frac\nu{\rho_*}-d\right)\norm{\nabla\Div h_m}_2^2
 +\frac{2\mu}{\rho_*}
          \norm{\nabla\times(\nabla\times h_m)}_2^2\right\},
\]
where the definition of $d$ also gives
$\nu/\rho_*-d\ge\nu/(2\rho_*)>0$.
The identity
\mbox{$\norm{D^2h_m}_2^2=
\norm{\nabla\Div h_m}_2^2+
\norm{\nabla\times(\nabla\times h_m)}_2^2$}
then shows that these terms control $D_0$.

We interpret the terms containing three derivatives of $r$
in the distributional sense by moving two derivatives onto
$h_m$:
\begin{equation}\label{g2:density-transfer}
\begin{aligned}
&\langle-d\mathcal L(\widetilde v\nabla r)
                          +d^2\nabla\Delta r,h_m\rangle\\
&\quad=-d\int\widetilde v\nabla r\cdot\mathcal Lh_m\dd y
                 +d^2\int\nabla r\cdot\Delta h_m\dd y.
\end{aligned}
\end{equation}
The resulting pairing is bounded by
$C\norm{\nabla r}_2D_0^{1/2}$. For the heat contribution,
\[
 \left\langle\kappa d\Delta\left(
       \frac{\widetilde u_1}{c_v\widetilde\rho}\partial_1r\right),
       h_E\right\rangle
 =\frac{\kappa d}{c_v}\int
       \frac{\widetilde u_1}{\widetilde\rho}
                     \partial_1r\,\Delta h_E\dd y,
\]
so the heat pairing is bounded by
$C\delta\norm{\nabla r}_2D_0^{1/2}$. For the viscous work
terms, apply \eqref{g2:viscous-work} with
\[
 w=\frac{g_m-\widetilde u_1re_1}{\widetilde\rho}
 \quad\hbox{and}\quad
 w=-d\widetilde\rho^{-1}\nabla r.
\]
Their $L^2$ norms are bounded, respectively, by
$C\norm g_2$ and $C\norm{\nabla r}_2$.
Terms containing derivatives of the profile coefficients are
estimated as in Steps~3--4 of \cite{HuangXu2026}. Inserting
$\norm g_2\le C(D_0^{1/2}+\delta D_s^{1/2})$ and
$\norm{\widetilde u_1'\nabla h}_2\le C\delta D_s^{1/2}$
into these estimates leads, for any $\eta>0$, to
\[
 \langle\mathcal B g,h\rangle
 \le-(2c_0-C\delta-\eta)D_0
       +C_\eta\delta^2D_s+C_\eta\norm{\nabla r}_2^2.
\]
Choose $\eta=c_0/2$, then take $\delta$ small enough that
$C\delta\le c_0/2$ and $C_\eta\delta^2\le c_E/2$.
This derivation of the diffusion estimate requires only
$\nabla r\in L^2$ and $h\in\dot H^1\cap\dot H^2$.

Substituting \eqref{g2:recovery-identities} into the definitions of
$\widehat u,\widehat\theta$ gives
\[
 \widehat u=\frac\rho{\widetilde\rho}\psi,
 \qquad
 \widehat\theta=\frac\rho{\widetilde\rho}\zeta
                  +\frac\rho{2c_v\widetilde\rho}|\psi|^2.
\]
The pressure and convective momentum fluxes decompose as follows:
\[
\begin{aligned}
 p-\widetilde\rho\bar\theta
 &=\bar\theta r+\rho\zeta
  =\bar\theta r+\widetilde\rho\widehat\theta
                                  -\frac\rho{2c_v}|\psi|^2,\\
 \rho u\otimes u-\widetilde\rho\widetilde u\otimes\widetilde u
 &=\widetilde u\otimes f_m+f_m\otimes\widetilde u
       -r\widetilde u\otimes\widetilde u+\rho\psi\otimes\psi.
\end{aligned}
\]
The scalar factor in the total energy flux satisfies
\[
\begin{aligned}
 &(\rho E+p)-(\widetilde\rho\widetilde E
                                      +\widetilde\rho\bar\theta)\\
 &\quad=\left((c_v+1)\bar\theta
                         +\frac12|\widetilde u|^2\right)r
       +\rho\left((c_v+1)\zeta+\widetilde u\cdot\psi
                                      +\frac12|\psi|^2\right).
\end{aligned}
\]
The energy component of the linearized flux is
\[
\begin{aligned}
 (DF_j(\widetilde U)f)_E
 &=\widetilde u_j
       (f_E+\bar\theta r+\widetilde\rho\widehat\theta)\\
 &\quad+\left((c_v+1)\bar\theta+\frac12|\widetilde u|^2\right)
                                (f_{m,j}-\widetilde u_jr).
\end{aligned}
\]
Subtracting the linearized quantities from their exact
counterparts gives
\begin{equation}\label{g2:quadratic-fluxes}
\begin{aligned}
Q_u&=\psi-\widehat u=-\frac r{\widetilde\rho}\psi,
 &Q_\theta&=\zeta-\widehat\theta
       =-\frac r{\widetilde\rho}\zeta
                    -\frac\rho{2c_v\widetilde\rho}|\psi|^2,\\
Q_m&=\rho\psi\otimes\psi-\frac\rho{2c_v}|\psi|^2I,\\
Q_E&=\rho\left\{(c_v+1)\zeta+\widetilde u\cdot\psi
                         +\tfrac12|\psi|^2\right\}\psi
                   -\frac\rho{2c_v}|\psi|^2\widetilde u.
\end{aligned}
\end{equation}
Thus the nonlinear part of the effective system is
\begin{equation}\label{g2:nonlinear-flux}
\begin{aligned}
(\mathcal N_g)_\rho={}&0,\\
(\mathcal N_g)_m={}&-\Div Q_m+\mathcal LQ_u,\\
(\mathcal N_g)_E={}&-\Div Q_E+\kappa\Delta Q_\theta\\
 &+\Div\{\mathcal T(Q_u)\widetilde u
              +\mathcal T(\widetilde u)Q_u+\mathcal T(\psi)\psi\}.
\end{aligned}
\end{equation}
On the fixed positive state neighborhood, \eqref{g2:quadratic-fluxes} implies
\[
 |Q_u|+|Q_\theta|+|Q_m|+|Q_E|\le C|f|^2.
\]
By interpolation, Sobolev embedding, and \eqref{g2:positive-equivalence},
\[
 \norm f_3\le\norm f_2^{1/2}\norm f_6^{1/2}
       \le C\norm f_{H^1}\le CN,\qquad
 \norm\psi_3\le C\norm\psi_{H^1}\le CN.
\]
H\"older's inequality and Sobolev embedding give
\begin{align}
\norm{(Q_u,Q_\theta,Q_m,Q_E)}_{6/5}
 &\le C\norm f_3\norm f_2\le CN\norm f_2,\nonumber\\
\norm{Q_u}_2+\norm{Q_\theta}_2
 &\le C\norm f_3\norm f_6\le CN\norm{\nabla f}_2.
\label{g2:flux-bounds}
\end{align}
Integrating the flux terms by parts transfers the derivatives
to $h$ and gives
\begin{align*}
|\langle-\Div(Q_m,Q_E),(h_m,h_E)\rangle|
 &\le\norm{(Q_m,Q_E)}_{6/5}\norm{\nabla(h_m,h_E)}_6
 \le CN\norm f_2D_0^{1/2},\\
|\langle\mathcal LQ_u,h_m\rangle|
       +|\langle\kappa\Delta Q_\theta,h_E\rangle|
 &\le C(\norm{Q_u}_2+\norm{Q_\theta}_2)D_0^{1/2}
 \le CN\norm{\nabla f}_2D_0^{1/2},\\
|\langle\Div\{\mathcal T(\psi)\psi\},h_E\rangle|
 &\le C\norm{\nabla\psi}_2\norm\psi_3\norm{\nabla h_E}_6
 \le CN\norm{\nabla\psi}_2D_0^{1/2}.
\end{align*}
For the viscous work terms, \cite[equation (3.14)]{HuangXu2026}
gives, in our notation,
\begin{equation}\label{g2:viscous-work}
\begin{aligned}
&\langle\Div\{\mathcal T(w)\widetilde u
                    +\mathcal T(\widetilde u)w\},h_E\rangle\\
&\quad=\int\widetilde u_1w\cdot
        \{(\mu+\lambda)\nabla\partial_1h_E
                                     +\mu e_1\Delta h_E\}\dd y\\
&\qquad+(\mu-\lambda)
                 \int\widetilde u_1'\sum_{j=2}^3w_j\partial_jh_E\dd y.
\end{aligned}
\end{equation}
The right-hand side is continuous for $w\in L^2$ and
$h_E\in\dot H^1\cap\dot H^2$, so the identity extends to these
spaces by density.
Set $w=Q_u$. The profile bounds
$\norm{\widetilde u_1}_\infty\le C\delta$ and
$\norm{\widetilde u_1'}_\infty\le C\delta^2$ give
\[
\begin{aligned}
&|\langle\Div\{\mathcal T(Q_u)\widetilde u
                       +\mathcal T(\widetilde u)Q_u\},h_E\rangle|\\
&\quad\le C\delta\norm{Q_u}_2D_0^{1/2}
             +C\delta^2\norm{Q_u}_{6/5}\norm{\nabla h_E}_6
 \\
&\quad\le CN\bigl(\delta\norm{\nabla f}_2
                   +\delta^2\norm f_2\bigr)D_0^{1/2}.
\end{aligned}
\]
Combining \eqref{g2:positive-equivalence} and \eqref{g2:recovery}
with the definition of $D_2$ gives
\[
\begin{aligned}
 D_0^{1/2}
       \bigl(\norm f_2+\norm{\nabla f}_2+\norm{\nabla\psi}_2\bigr)
 &\le CD_0^{1/2}(\norm Z_2+D_2^{1/2})\\
 &\le C(D_0+D_2+\norm Z_2^2).
\end{aligned}
\]
Summing the nonlinear pairing estimates, we obtain
\begin{equation}\label{g2:nonlinear-bound}
|\langle\mathcal N_g,h\rangle|
 \le CN D_0^{1/2}
       \bigl(\norm f_2+\norm{\nabla f}_2+\norm{\nabla\psi}_2\bigr)
 \le CN(D_0+D_2+\norm Z_2^2).
\end{equation}
Combining \eqref{g2:negative-identity}, \eqref{g2:effective-system},
\eqref{g2:Euler-pairing}, \eqref{g2:diffusive-pairing}, and
\eqref{g2:nonlinear-bound} yields
\[
 E_-'+c_0D_0+\frac{c_E}{2}D_s
 \le C\norm{\nabla r}_2^2
          +CN(D_0+D_2+\norm Z_2^2).
\]
Taking $c=\min\{c_0,c_E/2\}$ proves the lemma.
\end{proof}

\subsection{Second-order estimates}

\begin{lemma}\label{ft:H2-lemma}
If $\delta+N$ is sufficiently small, the parameter
$\sigma_2$ can be chosen so that $P_2\sim N^2$ and
\begin{equation}\label{ft:H2-energy}
P_2'+c_pD_2
 \le C_p(\delta+N)\bigl(D_2+\norm Z_2^2\bigr).
\end{equation}
\end{lemma}

\begin{proof}
From $v=\widetilde v e^{-\Phi}$, we compute
\[
\begin{aligned}
 \nabla(v-\widetilde v)
 &=-v\nabla\Phi
       +\widetilde v'(e^{-\Phi}-1)e_1,\\
 \nabla^2(v-\widetilde v)
 &=-v\nabla^2\Phi+v\nabla\Phi\otimes\nabla\Phi\\
 &\quad-\widetilde v'e^{-\Phi}
             (e_1\otimes\nabla\Phi+\nabla\Phi\otimes e_1)
       +\widetilde v''(e^{-\Phi}-1)e_1\otimes e_1.
\end{aligned}
\]
For $N\le1$,
\[
 \norm{v-\widetilde v}_{H^2}
 \le C\left(\norm\Phi_{H^2}+\norm{\nabla\Phi}_4^2\right)
 \le C(N+N^2)\le CN.
\]
The embeddings $H^2(\R^3)\hookrightarrow L^\infty$ and
$H^1(\R^3)\hookrightarrow L^6$ then give
\begin{equation}\label{ft:embeddings}
\begin{gathered}
\norm Z_\infty+\norm{\nabla Z}_3
 +\norm{v-\widetilde v}_{H^2}\le CN,\\
\norm{\nabla(\psi,\zeta)}_\infty
 +\norm{\nabla^2(\psi,\zeta)}_6
 +\norm{\nabla^2(\psi,\zeta)}_3\le C D_2^{1/2}.
\end{gathered}
\end{equation}
Throughout the estimates below, the profile coefficients are
bounded in $L^\infty$ by $C_\delta$, and their positive-order
derivatives are bounded by $C\delta^2$. In particular, they act
as bounded multipliers on the Sobolev spaces used below.

Set $(\Phi_\alpha,\psi_\alpha,\zeta_\alpha)
=\partial^\alpha(\Phi,\psi,\zeta)$ for $|\alpha|\le2$.
The differentiated equations in
\cite[equations (3.24)--(3.25)]{HuangXu2026}, restricted to
$|\alpha|\le2$, read
\begin{equation}\label{en:differentiated}
\left\{\begin{aligned}
&D_\tau\Phi_\alpha+\Div\psi_\alpha=F_{\Phi,\alpha},\\
&D_\tau\psi_\alpha+\theta\nabla\Phi_\alpha+\nabla\zeta_\alpha
                      -v\mathcal L\psi_\alpha=F_{\psi,\alpha},\\
&c_vD_\tau\zeta_\alpha+\theta\Div\psi_\alpha
                      -\kappa v\Delta\zeta_\alpha=F_{\zeta,\alpha},
\end{aligned}\right.
\end{equation}
The commutator terms are
\begin{equation}\label{en:commutators}
\begin{aligned}
F_{\Phi,\alpha}
 &=\partial^\alpha F_\Phi-[\partial^\alpha,u]\cdot\nabla\Phi,\\
F_{\psi,\alpha}
 &=\partial^\alpha F_\psi-[\partial^\alpha,u]\cdot\nabla\psi
   -[\partial^\alpha,\theta]\nabla\Phi
                          +[\partial^\alpha,v]\mathcal L\psi,\\
F_{\zeta,\alpha}
 &=\partial^\alpha F_\zeta-c_v[\partial^\alpha,u]\cdot\nabla\zeta
   -[\partial^\alpha,\theta]\Div\psi
                          +\kappa[\partial^\alpha,v]\Delta\zeta.
\end{aligned}
\end{equation}
Here and below,
$[\partial^\alpha,b]z=\partial^\alpha(bz)-b\partial^\alpha z$.

The energy identity in \cite[equations (3.26)--(3.28)]{HuangXu2026},
with $m=2$, uses the constant-state acoustic cancellation
\begin{equation}\label{en:acoustic-cancellation}
\begin{aligned}
&\theta_*\int_{\R^3}
       (\Phi_\alpha\Div\psi_\alpha
                           +\psi_\alpha\cdot\nabla\Phi_\alpha)\dd y\\
&\qquad+\int_{\R^3}
       (\psi_\alpha\cdot\nabla\zeta_\alpha
                             +\zeta_\alpha\Div\psi_\alpha)\dd y=0.
\end{aligned}
\end{equation}
Substituting the differentiated equations and using this
cancellation gives
\[
\begin{aligned}
&S_2'-\sum_{|\alpha|\le2}\left\{
       \int v\mathcal L\psi_\alpha\cdot\psi_\alpha
       +\frac\kappa{\theta_*}\int v\Delta\zeta_\alpha\,\zeta_\alpha
                                      \right\}\\
&\quad=\sum_{|\alpha|\le2}\biggl\{
 \frac12\int(\Div u)\left(
   \theta_*\Phi_\alpha^2+|\psi_\alpha|^2
                        +\frac{c_v}{\theta_*}\zeta_\alpha^2\right)\\
&\qquad+\int\left(\theta_*F_{\Phi,\alpha}\Phi_\alpha
       +F_{\psi,\alpha}\cdot\psi_\alpha
                       +\frac1{\theta_*}F_{\zeta,\alpha}\zeta_\alpha\right)\\
&\qquad-\int(\theta-\theta_*)\psi_\alpha\cdot\nabla\Phi_\alpha
       -\frac1{\theta_*}\int(\theta-\theta_*)
                               \zeta_\alpha\Div\psi_\alpha\biggr\}.
\end{aligned}
\]
For the density coupling, \cite[equation (3.30)]{HuangXu2026}
gives

\begin{equation}\label{ft:pressure-transfer}
-\int(\theta-\theta_*)\psi_\alpha\cdot\nabla\Phi_\alpha
 =\int(\theta-\theta_*)\Phi_\alpha\Div\psi_\alpha
  +\int\Phi_\alpha\psi_\alpha\cdot\nabla\theta.
\end{equation}
Since $\norm{\theta-\theta_*}_\infty\le C(\delta+N)$,
combining the first term on the right of \eqref{ft:pressure-transfer}
with the other temperature coupling term gives
\[
\begin{aligned}
&\left|\sum_{|\alpha|\le2}\int(\theta-\theta_*)
       \left(\Phi_\alpha-\frac{\zeta_\alpha}{\theta_*}\right)
                                      \Div\psi_\alpha\dd y\right|\\
&\quad\le C(\delta+N)\sum_{|\alpha|\le2}
       (\norm{\Phi_\alpha}_2+\norm{\zeta_\alpha}_2)
                                      \norm{\nabla\psi_\alpha}_2\\
&\quad\le C(\delta+N)(D_2+\norm Z_2^2).
\end{aligned}
\]
For $1\le|\alpha|\le2$, the perturbative part of the
last integral in \eqref{ft:pressure-transfer} is bounded by
\[
 \left|\int\Phi_\alpha\psi_\alpha\cdot\nabla\zeta\right|
 \le\norm{\Phi_\alpha}_2\norm{\psi_\alpha}_6
                   \norm{\nabla\zeta}_3\le CN D_2.
\]
The corresponding profile contribution is bounded by
$C\delta^2\norm{\Phi_\alpha}_2\norm{\psi_\alpha}_2$.
For $|\alpha|=0$, we have
\[
\begin{aligned}
 \left|\int\Phi\psi\cdot\nabla\zeta\dd y\right|
 &\le\norm\psi_\infty\norm\Phi_2\norm{\nabla\zeta}_2
       \le CN(D_2+\norm Z_2^2),\\
 \left|\int\bar\theta'\Phi\psi_1\dd y\right|
 &\le C\delta^2\norm\Phi_2\norm\psi_2
       \le C\delta^2\norm Z_2^2.
\end{aligned}
\]
The second-order transport commutator is
\[
 [\partial_i\partial_j,\psi]\cdot\nabla\Phi
   =(\partial_i\psi)\cdot\nabla\partial_j\Phi
       +(\partial_j\psi)\cdot\nabla\partial_i\Phi
       +(\partial_{ij}\psi)\cdot\nabla\Phi.
\]
For $1\le|\alpha|\le2$, the two terms with the highest
density derivatives obey
\begin{align*}
\left|\int(\Div\psi)\Phi_\alpha^2\right|
 &\le \norm{\nabla\psi}_\infty\norm{\Phi_\alpha}_2^2
 \le C D_2^{1/2}N D_2^{1/2}=CN D_2,\\
\left|\int\partial^\alpha\psi\cdot\nabla\Phi\,\Phi_\alpha\right|
 &\le\norm{\partial^\alpha\psi}_6
        \norm{\nabla\Phi}_3\norm{\Phi_\alpha}_2
 \le CN D_2.
\end{align*}
In the first estimate, one factor $\norm{\Phi_\alpha}_2$
is bounded by $N$ and the other by $D_2^{1/2}$. The term with one derivative on
$\psi$ and two on $\Phi$ has the same bound. For
$|\alpha|=0$, we use
$\norm Z_\infty\norm{\nabla\psi}_2\norm Z_2$.
For $z=\psi$ or $\zeta$ and $1\le|\alpha|\le2$,
\[
 \int|\Div\psi|\,|\partial^\alpha z|^2\dd y
       \le C\norm{\nabla\psi}_3\norm{\partial^\alpha z}_2\norm{\partial^\alpha z}_6
       \le CN D_2.
\]
The remaining second-order transport commutators are bounded by
\[
\begin{aligned}
 &\left|\left\langle
       [\partial_i\partial_j,\psi]\cdot\nabla z,\partial_{ij}z
                                                    \right\rangle\right|\\
 &\quad\le C\norm{\nabla\psi}_3
             \norm{\nabla^2z}_2\norm{\nabla^2z}_6
       +C\norm{\nabla^2\psi}_6
             \norm{\nabla z}_3\norm{\nabla^2z}_2
       \le CN D_2.
\end{aligned}
\]
The profile contribution is at most $C\delta^2D_2$, since
$\norm{\widetilde u_1'}_{W^{1,\infty}}\le C\delta^2$.
For the pressure commutator, expand
\[
 [\partial^\alpha,\theta]\nabla\Phi
   =\sum_{0<\beta\le\alpha}\binom\alpha\beta
       (\partial^\beta\theta)\nabla\partial^{\alpha-\beta}\Phi.
\]
When $|\alpha|=2$, Hölder's inequality gives
\[
\begin{aligned}
 \norm{\nabla\zeta}_3\norm{\nabla^2\Phi}_2
                         \norm{\nabla^2\psi}_6&\le CN D_2,\\
 \norm{\nabla^2\zeta}_6\norm{\nabla\Phi}_3
                         \norm{\nabla^2\psi}_2&\le CN D_2.
\end{aligned}
\]
The temperature acoustic commutator is treated in the same way:
\[
\begin{aligned}
 &\left|\int[\partial^\alpha,\zeta]\Div\psi\,\zeta_\alpha\dd y\right|\\
 &\quad\le C\norm{\nabla\zeta}_3
              \norm{\nabla^2\psi}_2\norm{\nabla^2\zeta}_6
       +C\norm{\nabla^2\zeta}_6
              \norm{\nabla\psi}_3\norm{\nabla^2\zeta}_2
 \le CN D_2.
\end{aligned}
\]

For $1\le|\alpha|\le2$, the commutators with the
profile temperature are
\[
\begin{aligned}
 [\partial^\alpha,\bar\theta]\nabla\Phi
 &=\sum_{j=1}^{\alpha_1}\binom{\alpha_1}{j}
       \bar\theta^{(j)}\nabla\partial^{\alpha-je_1}\Phi,\\
 [\partial^\alpha,\bar\theta]\Div\psi
 &=\sum_{j=1}^{\alpha_1}\binom{\alpha_1}{j}
       \bar\theta^{(j)}\Div\partial^{\alpha-je_1}\psi.
\end{aligned}
\]
The $L^\infty$ bounds on the profile derivatives give
\[
\begin{aligned}
 &\left|\left\langle
       [\partial^\alpha,\bar\theta]\nabla\Phi,\psi_\alpha
                                              \right\rangle\right|\\
 &\quad\le C\sum_{j=1}^{\alpha_1}
       \norm{\bar\theta^{(j)}}_\infty
       \norm{\nabla\partial^{\alpha-je_1}\Phi}_2
       \norm{\psi_\alpha}_2\\
 &\quad\le C\delta^2\norm{\nabla\Phi}_{H^1}
                            \norm{\nabla\psi}_{H^1}\le C\delta^2D_2,\\
 &\left|\left\langle
       [\partial^\alpha,\bar\theta]\Div\psi,\zeta_\alpha
                                              \right\rangle\right|\\
 &\quad\le C\sum_{j=1}^{\alpha_1}
       \norm{\bar\theta^{(j)}}_\infty
       \norm{\nabla\partial^{\alpha-je_1}\psi}_2
       \norm{\zeta_\alpha}_2\\
 &\quad\le C\delta^2\norm{\nabla\psi}_{H^1}
                            \norm{\nabla\zeta}_{H^1}\le C\delta^2D_2.
\end{aligned}
\]
Both commutators vanish at order zero.

For each component $z$ of $\psi$ and for $z=\zeta$, the
diffusion commutator expands as
\[
\begin{aligned}
 [\partial_i\partial_j,v]\partial_{\ell m}z
 &=(\partial_i v)\partial_{j\ell m}z
       +(\partial_j v)\partial_{i\ell m}z
       +(\partial_{ij}v)\partial_{\ell m}z.
\end{aligned}
\]
Decompose $v=\widetilde v+(v-\widetilde v)$. The terms
containing $v-\widetilde v$ satisfy
\begin{align*}
\left|\int\nabla(v-\widetilde v)\,\nabla^3z\,\nabla^2z\right|
 &\le\norm{\nabla(v-\widetilde v)}_3
        \norm{\nabla^3z}_2\norm{\nabla^2z}_6\le CN D_2,\\
\left|\int\nabla^2(v-\widetilde v)\,\nabla^2z\,\nabla^2z\right|
 &\le\norm{\nabla^2(v-\widetilde v)}_2
        \norm{\nabla^2z}_3\norm{\nabla^2z}_6\le CN D_2.
\end{align*}
The remaining profile terms satisfy
\[
\begin{aligned}
 &\int|\widetilde v'|\,|\nabla^3z|\,|\nabla^2z|\dd y
       +\int|\widetilde v''|\,|\nabla^2z|^2\dd y\\
 &\qquad\le C\delta^2\left(
       \norm{\nabla^3z}_2\norm{\nabla^2z}_2
                                    +\norm{\nabla^2z}_2^2\right)
 \le C\delta^2D_2.
\end{aligned}
\]
For the principal diffusion terms, integration by parts gives
\[
\begin{aligned}
-\int v\mathcal L\psi_\alpha\cdot\psi_\alpha
 &=\mu\int v|\nabla\psi_\alpha|^2
           +(\mu+\lambda)\int v|\Div\psi_\alpha|^2\\
 &\quad+\mu\sum_{j=1}^3\int(\partial_jv)
                         \partial_j\psi_\alpha\cdot\psi_\alpha
       +(\mu+\lambda)\int(\nabla v\cdot\psi_\alpha)\Div\psi_\alpha,\\
-\int v\Delta\zeta_\alpha\,\zeta_\alpha
 &=\int v|\nabla\zeta_\alpha|^2
                      +\int\zeta_\alpha\nabla v\cdot\nabla\zeta_\alpha.
\end{aligned}
\]
Since $v$ has a positive lower bound and
$\mu+\lambda\ge\mu/3>0$, the principal terms control
the velocity and temperature dissipation:
\[
\begin{aligned}
&\sum_{|\alpha|\le2}\left\{
 \mu\int v|\nabla\psi_\alpha|^2
       +(\mu+\lambda)\int v|\Div\psi_\alpha|^2
       +\frac\kappa{\theta_*}\int v|\nabla\zeta_\alpha|^2\right\}\\
&\quad\ge c\norm{\nabla(\psi,\zeta)}_{H^2}^2.
\end{aligned}
\]
For either $z=\psi$ or $z=\zeta$, the coefficient derivative
terms obey
\[
 \int|\nabla(v-\widetilde v)|\,|\partial^\alpha z|\,|\nabla \partial^\alpha z|
 \le\norm{\nabla(v-\widetilde v)}_3\norm{\partial^\alpha z}_6
             \norm{\nabla \partial^\alpha z}_2
 \le CN\norm{\nabla \partial^\alpha z}_2^2.
\]
The profile terms are bounded using
$\norm{\widetilde v'}_\infty\le C\delta^2$. Their sum over
$|\alpha|\le2$ is at most
$C\delta^2(D_2+\norm Z_2^2)$.

For the viscous heating term, linearity of $\mathcal T$ gives
\[
\begin{aligned}
 \partial_j\{v\mathcal T(\psi):\nabla\psi\}
 &=(\partial_jv)\mathcal T(\psi):\nabla\psi\\
 &\quad+v\mathcal T(\partial_j\psi):\nabla\psi
                   +v\mathcal T(\psi):\nabla\partial_j\psi.
\end{aligned}
\]
Decomposing $v=\widetilde v+(v-\widetilde v)$ and applying
H\"older's inequality gives
\begin{align}
&\norm{\nabla\{v\mathcal T(\psi):\nabla\psi\}}_2\nonumber\\
&\quad\le C\norm{\nabla(v-\widetilde v)}_6\norm{\nabla\psi}_6^2
 +C\norm v_\infty\norm{\nabla^2\psi}_3\norm{\nabla\psi}_6
 +C\norm{\widetilde v'}_\infty\norm{\nabla\psi}_4^2
 \le CN D_2^{1/2}.\label{ft:heating}
\end{align}
For the last term, use
$\norm{\nabla\psi}_4^2\le
\norm{\nabla\psi}_3\norm{\nabla\psi}_6$. The three terms
in \eqref{ft:heating} are therefore bounded as follows:
\[
\begin{aligned}
\norm{\nabla(v-\widetilde v)}_6\norm{\nabla\psi}_6^2
 &\le CN^2D_2^{1/2},\\
\norm{\nabla^2\psi}_3\norm{\nabla\psi}_6
 &\le CN D_2^{1/2},\\
\norm{\widetilde v'}_\infty\norm{\nabla\psi}_4^2
 &\le C\delta^2N D_2^{1/2}.
\end{aligned}
\]
Since $N\le1$, this proves \eqref{ft:heating} with only first
derivatives of $v$.
For $|\alpha|=2$, write
$\alpha=\beta+e_j$, where $|\beta|=1$.
Moving one derivative onto the temperature perturbation gives
\begin{align*}
&\left|\int\partial^\alpha\{v\mathcal T(\psi):\nabla\psi\}
                                      \zeta_\alpha\dd y\right|\\
&\quad=\left|\int\partial^\beta\{v\mathcal T(\psi):\nabla\psi\}
                              \partial_j\zeta_\alpha\dd y\right|
 \le C\norm{\nabla\{v\mathcal T(\psi):\nabla\psi\}}_2
                     \norm{\nabla^3\zeta}_2
 \le CN D_2.
\end{align*}
For $|\alpha|=1$, pair \eqref{ft:heating} with
$\nabla\zeta$; for $|\alpha|=0$, the contribution is bounded
by $C\norm\zeta_\infty\norm{\nabla\psi}_2^2$.
The mixed viscous heating term involving the profile is
\[
 2v\mathcal T(\widetilde u):\nabla\psi
 =2v\widetilde u_1'
       (\nu\partial_1\psi_1+\lambda\partial_2\psi_2
                                             +\lambda\partial_3\psi_3).
\]
For $1\le i,j,\ell\le3$, direct differentiation gives
\[
\begin{aligned}
 \partial_{ij}(v\partial_\ell\psi_\ell)
 &=v\partial_{ij\ell}\psi_\ell
       +(\partial_i v)\partial_{j\ell}\psi_\ell
       +(\partial_j v)\partial_{i\ell}\psi_\ell
       +(\partial_{ij}v)\partial_\ell\psi_\ell.
\end{aligned}
\]
After writing $v=\widetilde v+(v-\widetilde v)$, the terms
containing the coefficient perturbation obey
\[
\begin{aligned}
 &\norm{(v-\widetilde v)\nabla^3\psi}_2
       +\norm{\nabla(v-\widetilde v)\,\nabla^2\psi}_2
       +\norm{\nabla^2(v-\widetilde v)\,\nabla\psi}_2\\
 &\quad\le\norm{v-\widetilde v}_\infty\norm{\nabla^3\psi}_2
       +\norm{\nabla(v-\widetilde v)}_3\norm{\nabla^2\psi}_6
       +\norm{\nabla^2(v-\widetilde v)}_2\norm{\nabla\psi}_\infty\\
 &\quad\le CN D_2^{1/2}.
\end{aligned}
\]
The profile terms are bounded by
\[
\begin{aligned}
 &\norm{\widetilde v\nabla^3\psi}_2
       +\norm{\widetilde v'\nabla^2\psi}_2
       +\norm{\widetilde v''\nabla\psi}_2\\
 &\quad\le C\norm{\nabla\psi}_{H^2}.
\end{aligned}
\]
At orders zero and one, the required bounds are
\[
\begin{aligned}
 \norm{v\nabla\psi}_2
 &\le\norm v_\infty\norm{\nabla\psi}_2,\\
 \norm{\nabla(v\nabla\psi)}_2
 &\le\norm v_\infty\norm{\nabla^2\psi}_2
       +\norm{\widetilde v'}_\infty\norm{\nabla\psi}_2
       +\norm{\nabla(v-\widetilde v)}_3\norm{\nabla\psi}_6.
\end{aligned}
\]
Differentiating the factor $\widetilde u_1'$ gives
\[
\begin{aligned}
 \partial_{ij}(\widetilde u_1'v\partial_\ell\psi_\ell)
 &=\widetilde u_1'\partial_{ij}(v\partial_\ell\psi_\ell)\\
 &\quad+\widetilde u_1''\bigl\{
       \delta_{i1}\partial_j(v\partial_\ell\psi_\ell)
       +\delta_{j1}\partial_i(v\partial_\ell\psi_\ell)\bigr\}\\
 &\quad+\delta_{i1}\delta_{j1}
                    \widetilde u_1'''v\partial_\ell\psi_\ell.
\end{aligned}
\]
The bound
$\norm{\widetilde u_1'}_{W^{2,\infty}}\le C\delta^2$, together
with $N\le1$, now gives
\[
 \norm{2v\mathcal T(\widetilde u):\nabla\psi}_{H^2}
       \le C\delta^2\norm{\nabla\psi}_{H^2}
       \le C\delta^2D_2^{1/2}.
\]
Pairing this term with $\zeta$ in $H^2$ gives
$C\delta^2N D_2^{1/2}\le C\delta^2(D_2+\norm Z_2^2)$,
where we used \mbox{$N^2\le C(D_2+\norm Z_2^2)$}.
The remaining profile source terms satisfy
\[
\begin{aligned}
 \norm{F_\Phi}_{H^2}+\norm{F_\psi}_{H^2}
       &\le C\delta^2N,\\
 \norm{F_\zeta-v\mathcal T(\psi):\nabla\psi
                   -2v\mathcal T(\widetilde u):\nabla\psi}_{H^2}
       &\le C\delta^2N.
\end{aligned}
\]
For the density source, we have
\[
\begin{aligned}
 \partial_{ij}\left(\frac{\widetilde\rho'}{\widetilde\rho}\psi_1\right)
 &=\frac{\widetilde\rho'}{\widetilde\rho}\partial_{ij}\psi_1\\
 &\quad+\left(\frac{\widetilde\rho'}{\widetilde\rho}\right)'
       (\delta_{i1}\partial_j\psi_1+\delta_{j1}\partial_i\psi_1)\\
 &\quad+\delta_{i1}\delta_{j1}
       \left(\frac{\widetilde\rho'}{\widetilde\rho}\right)''\psi_1.
\end{aligned}
\]
The positive lower bound for the density, together with the
profile derivative estimates, gives
\[
\begin{aligned}
 \norm{\frac{\widetilde\rho'}{\widetilde\rho}\psi_1}_{H^2}
 &\le C\sum_{j=0}^2
       \norm{\partial_1^j(\widetilde\rho'/\widetilde\rho)}_\infty
       \norm{\psi_1}_{H^{2-j}}\\
 &\le C\delta^2N.
\end{aligned}
\]
For the heat source involving $\bar\theta''$,
\[
\begin{aligned}
 \partial_{ij}((v-\widetilde v)\bar\theta'')
 &=\bar\theta''\partial_{ij}(v-\widetilde v)
       +\bar\theta'''(\delta_{i1}\partial_j(v-\widetilde v)
                            +\delta_{j1}\partial_i(v-\widetilde v))\\
 &\quad+\delta_{i1}\delta_{j1}\bar\theta^{(4)}(v-\widetilde v).
\end{aligned}
\]
It follows that
\[
\begin{aligned}
 \norm{(v-\widetilde v)\bar\theta''}_{H^2}
 &\le C\norm{\bar\theta''}_\infty\norm{v-\widetilde v}_{H^2}
       +C\norm{\bar\theta'''}_\infty\norm{v-\widetilde v}_{H^1}\\
 &\quad+C\norm{\bar\theta^{(4)}}_\infty\norm{v-\widetilde v}_2
 \le C\delta^2N.
\end{aligned}
\]
The coefficients $\widetilde u_1'$, $\widetilde u_1''$ and
$\bar\theta'$ have $W^{2,\infty}$ norms bounded by
$C\delta^2$. Also,
\[
 \norm{(\widetilde u_1')^2}_{W^{2,\infty}}\le C\delta^4.
\]
These coefficient bounds give
\[
\begin{aligned}
 &\norm{\widetilde u_1'\psi_1}_{H^2}
       +\norm{\frac{\widetilde\rho'}{\widetilde\rho}\zeta}_{H^2}
       +\norm{(v-\widetilde v)\widetilde u_1''}_{H^2}\\
 &\quad+\norm{\bar\theta'\psi_1}_{H^2}
       +\norm{\widetilde u_1'\zeta}_{H^2}
       +\norm{(v-\widetilde v)(\widetilde u_1')^2}_{H^2}
 \le C\delta^2N.
\end{aligned}
\]
For the heat source, pairing with $\zeta$ gives
\[
\begin{aligned}
 &\sum_{|\alpha|\le2}\left|
       \int\partial^\alpha((v-\widetilde v)\bar\theta'')
                                      \zeta_\alpha\dd y\right|\\
 &\quad\le\norm{(v-\widetilde v)\bar\theta''}_{H^2}\norm\zeta_{H^2}
       \le C\delta^2N^2
       \le C\delta^2(D_2+\norm Z_2^2).
\end{aligned}
\]
The remaining source terms are estimated in the same way, using
the $H^2$ bounds above. Summing over $|\alpha|\le2$, we obtain
\begin{equation}\label{ft:symmetric-H2}
S_2'+c_0\norm{\nabla(\psi,\zeta)}_{H^2}^2
 \le C(\delta+N)(D_2+\norm Z_2^2).
\end{equation}

For the cross term, \cite[equation (3.31)]{HuangXu2026},
with $m=2$, gives
\begin{equation}\label{g2:interactive-identity}
\begin{aligned}
&I_2'+\sum_{|\alpha|\le1}\int\theta|\nabla\Phi_\alpha|^2\\
&\quad=\sum_{|\alpha|\le1}\left\{
 \norm{\Div\psi_\alpha}_2^2
 -\int\nabla\zeta_\alpha\cdot\nabla\Phi_\alpha
 +\int v\mathcal L\psi_\alpha\cdot\nabla\Phi_\alpha\right\}\\
&\qquad+\sum_{|\alpha|\le1}\int_{\R^3}\biggl\{
       (\Div u)\psi_\alpha\cdot\nabla\Phi_\alpha
       -\sum_{j,k=1}^3\psi_{\alpha,j}(\partial_j u_k)
                                           \partial_k\Phi_\alpha\\
&\hspace{42mm}+F_{\psi,\alpha}\cdot\nabla\Phi_\alpha
                         -F_{\Phi,\alpha}\Div\psi_\alpha\biggr\}\dd y.
\end{aligned}
\end{equation}
By \eqref{en:commutators}, for $|\alpha|\le1$,
\[
\norm{F_{\Phi,\alpha}}_2+\norm{F_{\psi,\alpha}}_2
 \le C(\delta^2+N)\bigl(D_2^{1/2}+\norm Z_2\bigr).
\]
For the diffusion and pressure commutators, H\"older's
inequality in $L^3\times L^6$ gives
\[
\begin{aligned}
\norm{\nabla(v-\widetilde v)\,\nabla^2\psi}_2
 &\le\norm{\nabla(v-\widetilde v)}_3\norm{\nabla^2\psi}_6
                                                  \le CN D_2^{1/2},\\
\norm{\nabla\zeta\,\nabla\Phi}_2
 &\le\norm{\nabla\zeta}_3\norm{\nabla\Phi}_6
                                                  \le CN D_2^{1/2}.
\end{aligned}
\]
The purely perturbative transport terms satisfy
\[
\begin{aligned}
 \sum_{j=1}^3\bigl(
       \norm{(\partial_j\psi)\cdot\nabla\Phi}_2
       +\norm{(\partial_j\psi)\cdot\nabla\psi}_2\bigr)
 &\le C\norm{\nabla\psi}_3\norm{\nabla(\Phi,\psi)}_6\\
 &\le CN D_2^{1/2}.
\end{aligned}
\]
The transport remainder contains products of
$\partial^\alpha\psi$, $\nabla u$, and
$\nabla\Phi_\alpha$. Its perturbative part satisfies
\[
\begin{aligned}
 \sum_{|\alpha|=1}\int
       |\psi_\alpha|\,|\nabla\psi|\,|\nabla\Phi_\alpha|\dd y
 &\le C\norm{\nabla\psi}_2\norm{\nabla\psi}_\infty
                                      \norm{\nabla^2\Phi}_2
       \le CN D_2,\\
 \int|\psi|\,|\nabla\psi|\,|\nabla\Phi|\dd y
 &\le\norm\psi_\infty\norm{\nabla\psi}_2\norm{\nabla\Phi}_2
       \le CN D_2.
\end{aligned}
\]
The profile terms are bounded by
$C\delta^2(D_2+\norm Z_2^2)$. Together with the preceding
$L^2$ estimates, this gives
\[
\begin{aligned}
&\sum_{|\alpha|\le1}\left(
       \norm{F_{\psi,\alpha}}_2\norm{\nabla\Phi_\alpha}_2
       +\norm{F_{\Phi,\alpha}}_2\norm{\Div\psi_\alpha}_2\right)\\
&\quad\le C(\delta^2+N)
       (D_2^{1/2}+\norm Z_2)D_2^{1/2}
       \le C(\delta+N)(D_2+\norm Z_2^2).
\end{aligned}
\]
The transport and source estimates above imply
\[
\begin{aligned}
&\Biggl|\sum_{|\alpha|\le1}\int_{\R^3}\biggl\{
       (\Div u)\psi_\alpha\cdot\nabla\Phi_\alpha
       -\sum_{j,k=1}^3\psi_{\alpha,j}(\partial_j u_k)
                                           \partial_k\Phi_\alpha\\
&\hspace{35mm}+F_{\psi,\alpha}\cdot\nabla\Phi_\alpha
                         -F_{\Phi,\alpha}\Div\psi_\alpha\biggr\}\dd y\Biggr|\\
&\quad\le C(\delta+N)(D_2+\norm Z_2^2).
\end{aligned}
\]
We estimate the two remaining mixed terms by Young's
inequality:
\[
\begin{aligned}
&\sum_{|\alpha|\le1}\biggl(
 \left|\int\nabla\zeta_\alpha\cdot\nabla\Phi_\alpha\right|
 +\left|\int v\mathcal L\psi_\alpha\cdot\nabla\Phi_\alpha\right|
                                                        \biggr)\\
&\qquad\le\eta\norm{\nabla\Phi}_{H^1}^2
       +C_\eta\norm{\nabla(\psi,\zeta)}_{H^2}^2.
\end{aligned}
\]
For sufficiently small $\delta+N$, we have
$\theta\ge\theta_*/2$. Taking $\eta=\theta_*/4$ in
\eqref{g2:interactive-identity} gives
\[
I_2'+c_1\norm{\nabla\Phi}_{H^1}^2
 \le C_1\norm{\nabla(\psi,\zeta)}_{H^2}^2
       +C(\delta+N)(D_2+\norm Z_2^2).
\]
Cauchy--Schwarz bounds the cross term by
\[
 |I_2|\le\sum_{|\alpha|\le1}
       \norm{\psi_\alpha}_2\norm{\nabla\Phi_\alpha}_2
       \le C S_2,\qquad S_2\sim N^2.
\]
Choose $\sigma_2>0$ small enough that
$\sigma_2C\le1/2$ and $\sigma_2C_1\le c_0/2$.
Then
\[
 \frac12S_2\le P_2=S_2+\sigma_2I_2\le\frac32S_2.
\]
Adding $\sigma_2$ times the estimate for $I_2'$ to
\eqref{ft:symmetric-H2} gives
\[
\begin{aligned}
&P_2'+\frac{c_0}{2}\norm{\nabla(\psi,\zeta)}_{H^2}^2
       +\sigma_2c_1\norm{\nabla\Phi}_{H^1}^2\\
&\qquad\le C(\delta+N)(D_2+\norm Z_2^2).
\end{aligned}
\]
Thus \eqref{ft:H2-energy} holds with
$c_p=\min\{c_0/2,\sigma_2c_1\}>0$. The energy pairings are
justified by the Fourier approximation in Lemma~\ref{g2:local}.
Both the pressure and viscous heating estimates use at most two
derivatives of the density. Third derivatives of velocity and
temperature are controlled by the time-integrated
dissipation.
\end{proof}

\subsection{Proof of the a priori estimate}

\begin{proof}[Proof of Proposition~\ref{g2:global-unit}]
We use the exponent $k$ in the elliptic energy and the
coefficient $\sigma_2$ chosen in Lemma~\ref{ft:H2-lemma}.
For $\delta\le1$, the recovery estimates
\eqref{g2:recovery} give
\[
\begin{aligned}
 \norm{\nabla r}_2^2
 &\le CD_2+C\delta^4(D_0+\delta^2D_s),\\
 D_0+D_2+\norm Z_2^2
 &\le C(D_0+D_2+\delta^2D_s)
       \le C(D_0+D_s+D_2).
\end{aligned}
\]
Inserting these bounds into
\eqref{g2:negative-inequality}, we obtain
\[
\begin{aligned}
 E_-'+c(D_0+D_s)
 &\le C_-D_2+C\delta^4(D_0+\delta^2D_s)\\
 &\quad+CN(D_0+D_s+D_2).
\end{aligned}
\]
We reduce the shock strength, if necessary, so that
$C\delta^4\le c/2$. Then
\[
 C\delta^4(D_0+\delta^2D_s)
       \le C\delta^4(D_0+D_s)
       \le\frac c2(D_0+D_s).
\]
This term is absorbed on the left, leaving
\[
 E_-'+c_-(D_0+D_s)
       \le C_-D_2+CN(D_0+D_s+D_2).
\]
Here $c_->0$ is fixed. The recovery estimates also give
\[
 D_2+\norm Z_2^2
       \le C(D_0+D_2+\delta^2D_s)
       \le C(D_0+D_s+D_2).
\]
Using this bound in \eqref{ft:H2-energy} and enlarging $C_p$
if necessary, we have
\[
 P_2'+c_pD_2
       \le C_p(\delta+N)(D_0+D_s+D_2).
\]
Set
\[
 L=\max\left\{1,\frac{C_-+1}{c_p}\right\}.
\]
Since $Lc_p-C_-\ge1$, we multiply the inequality for $P_2'$
by $L$ and add the negative-order estimate to obtain
\[
\begin{aligned}
 &(E_-+LP_2)'+c_-(D_0+D_s)+(Lc_p-C_-)D_2\\
 &\qquad\le\{CN+LC_p(\delta+N)\}(D_0+D_s+D_2).
\end{aligned}
\]
Keeping $L$ fixed, take $\delta_*>0$ and $N_*>0$ small
enough that all preceding estimates hold and
\[
 LC_p\delta_*\le\frac14\min\{c_-,1\},\qquad
 (C+LC_p)N_*\le\frac14\min\{c_-,1\}.
\]
Decrease $N_*$ further, if necessary, to lie below the continuation
radius in Lemma~\ref{g2:local}. Then $0<\delta\le\delta_*$ and
$N\le N_*$ imply
\[
 CN+LC_p(\delta+N)
       =LC_p\delta+(C+LC_p)N
       \le CN_*+LC_p(\delta_*+N_*)
       \le\frac12\min\{c_-,1\}.
\]
The right-hand side can therefore be absorbed, giving
\begin{equation}\label{g2:closed-energy}
\mathcal E:=E_-+LP_2,\qquad
\mathcal E'+c(D_0+D_s+D_2)\le0.
\end{equation}
Here $c=\frac12\min\{c_-,1\}>0$, and all constants are
independent of $S$.

Since $g-f=(0,d\nabla r,0)^{\mathsf T}$ and
$\norm{\nabla r}_{\dot H^{-1}}\le C\norm r_2\le C\norm Z_2$, we have
\[
\begin{aligned}
 \norm g_{\dot H^{-1}}
     &\le\norm f_{\dot H^{-1}}+C\norm Z_2,\\
 \norm f_{\dot H^{-1}}
     &\le\norm g_{\dot H^{-1}}+C\norm Z_2.
\end{aligned}
\]
Since $L$ is fixed, these bounds and the equivalence
\eqref{g2:negative-energy}, together with
$P_2\sim\norm Z_{H^2}^2$, yield
\[
\begin{aligned}
 \norm f_{\dot H^{-1}}^2+\norm Z_{H^2}^2
 &\le C\bigl(\norm g_{\dot H^{-1}}^2+\norm Z_{H^2}^2\bigr)
       \le C\mathcal E,\\
 \mathcal E
 &\le C\bigl(\norm g_{\dot H^{-1}}^2+\norm Z_{H^2}^2\bigr)
       \le C\bigl(\norm f_{\dot H^{-1}}^2+\norm Z_{H^2}^2\bigr).
\end{aligned}
\]
At the initial time, the Sobolev norm equivalence
\eqref{g2:positive-equivalence} gives
\[
\begin{aligned}
 \mathcal E(0)
 &\le C\bigl(\norm{g_0}_{\dot H^{-1}}^2
                            +\norm{Z_0}_{H^2}^2\bigr)\\
 &\le C\bigl(\norm{f_0}_{\dot H^{-1}}^2
                     +\norm{r_0}_2^2+\norm{Z_0}_{H^2}^2\bigr)\\
 &\le C_\delta\bigl(\norm{f_0}_{\dot H^{-1}}^2
                            +\norm{f_0}_{H^2}^2\bigr).
\end{aligned}
\]
Consequently,
\begin{equation}\label{g2:energy-equivalence}
 \mathcal E\sim\norm f_{\dot H^{-1}}^2+\norm Z_{H^2}^2,
 \qquad
 \mathcal E(0)\le C_\delta
       \bigl(\norm{f_0}_{\dot H^{-1}}^2+\norm{f_0}_{H^2}^2\bigr).
\end{equation}
Integrating \eqref{g2:closed-energy} gives
\[
 \mathcal E(S)+c\int_0^S(D_0+D_s+D_2)\dd\tau\le\mathcal E(0).
\]
Applying the same estimate on $[0,\tau]$ for
$0\le\tau\le S$ gives
$\sup_{0\le\tau\le S}\mathcal E(\tau)\le\mathcal E(0)$.
Since $\mathcal E(S)\ge0$, the integrated dissipation also
satisfies
\[
 \int_0^S(D_0+D_s+D_2)\dd\tau
       \le c^{-1}\mathcal E(0).
\]
Finally, \eqref{g2:recovery} and the definition of the
second-order dissipation in \eqref{ft:P2} give
\[
\begin{aligned}
 \norm Z_{H^2}^2+\norm{\nabla(\psi,\zeta)}_{H^2}^2
 &\le C(D_0+D_2+\delta^2D_s)+D_2\\
 &\le C(D_0+D_s+D_2).
\end{aligned}
\]
Combining these bounds with \eqref{g2:energy-equivalence}
gives \eqref{g2:unit-energy} with a constant independent of $S$:
\begin{align*}
&\sup_{0\le\tau\le S}
       \bigl(\norm{f(\tau)}_{\dot H^{-1}}^2+\norm{Z(\tau)}_{H^2}^2\bigr)\\
&\quad+\int_0^S\left\{
       \norm Z_{H^2}^2+\norm{\nabla(\psi,\zeta)}_{H^2}^2\right\}\dd\tau\\
&\quad\le C\sup_{0\le\tau\le S}\mathcal E(\tau)
       +C\int_0^S(D_0+D_s+D_2)\dd\tau\\
&\quad\le C(1+c^{-1})\mathcal E(0)\\
&\quad\le C_\delta\bigl(
       \norm{f_0}_{\dot H^{-1}}^2+\norm{f_0}_{H^2}^2\bigr).\tag*{\qedhere}
\end{align*}
\end{proof}

\setlength{\bibsep}{0.5em}

\end{document}